\documentclass[preprint,3p,12pt]{elsarticle}
\usepackage[utf8]{inputenc}
\usepackage{textcomp}

\usepackage{graphicx}
\usepackage{amsmath,amssymb}
\usepackage[version=4]{mhchem}
\usepackage{siunitx}
\usepackage{longtable,tabularx}
\usepackage[symbol]{footmisc}

\usepackage{bm}
\usepackage{amsfonts}
\usepackage{amscd}
\usepackage{hyperref}
\usepackage{pstricks}
\usepackage{subcaption}
\usepackage{framed,fancybox}
\usepackage{bbm}
\usepackage{multirow}
\usepackage{threeparttable}
\usepackage{rotating}
\usepackage{stfloats}
\usepackage{color,soul}
\usepackage{float}
\usepackage{afterpage}
\usepackage{listings}
\usepackage{xcolor}

\usepackage{amsthm}
\usepackage{xcolor}

\biboptions{square,numbers,sort&compress}
\newcommand{\m}[1]{{\bf{#1}}}

\newcommand{\dt}{\text{d}t}
\newcommand{\mb}[1]{\mathbf{#1}}

\definecolor{codegreen}{rgb}{0,0.6,0}
\definecolor{codegray}{rgb}{0.5,0.5,0.5}

\lstdefinestyle{mystyle}{
  commentstyle=\color{codegreen},
  numberstyle=\tiny\color{codegray},
  keywordstyle=\color{blue},
  basicstyle=\ttfamily\small,
  breakatwhitespace=false,         
  breaklines=true,                 
  keepspaces=true,                 
  numbers=left,                    
  numbersep=5pt,                  
  showspaces=false,                
  showstringspaces=false,
  showtabs=false,                  
  tabsize=2
}

\allowdisplaybreaks

\journal{Aerospace Science and Technology}

\begin{document}

\begin{frontmatter}

\title{{\bf Computational Method for Desensitized Optimal} \\{\bf Guidance Using Direct Collocation}\footnote[2]{Portions of this work were presented at the 2024 Conference on Decision and Control (CDC) in Milan, Italy (https://doi.org/10.1109/CDC56724.2024.10886014) while other portions were presented at the 2025 AAS/AIAA Spaceflight Mechanics Meeting in Lihue, Hawaii.}}

\author[ufaffil]{Katrina L.~Winkler\fnref{label1}}
\ead{k.winkler@ufl.edu}
\fntext[label1]{Ph.D. Candidate, Department of Mechanical and Aerospace Engineering.}
\affiliation[ufaffil]{organization={University of Florida},
            city={Gainesville},
            postcode={32611},
            state={FL},
            country={USA}}

\author[ufaffil]{Anil V. Rao\corref{cor1}\fnref{label2}}
\ead{anilvrao@ufl.edu}
\cortext[cor1]{Corresponding author.}
\fntext[label2]{Professor, Department of Mechanical and Aerospace Engineering. AIAA Associate Fellow. AAS Fellow.}

\begin{abstract}
A computational method is developed for desensitized optimal guidance using adaptive Gaussian quadrature collocation. The method computes a reference trajectory that reduces the sensitivity to uncertainties in the dynamic model by augmenting the objective functional to explicitly penalize the sensitivity of the state with respect to uncertain parameters. Using this desensitized reference trajectory as a starting point, the desensitized optimal guidance method developed in this paper computes a new optimal control on the remaining horizon at specified guidance update times. This shrinking horizon optimal control problem is solved using a Legendre-Gauss-Radau collocation method where, at each guidance update, a reduced-horizon mesh is determined by remapping the mesh to the remaining horizon and deleting the portion of the mesh associated with the expired portion of the horizon.  The resulting guidance solution is found to improve robustness to external disturbances and modeling errors.  The method is demonstrated on two numerical examples. The first example is Zermelo's navigation problem which illustrates the behavior of the method on a simple example. The second example is an atmospheric reentry problem that demonstrates the performance of the method on a more complex problem. For both examples, the dynamics are simulated in the presence of parameter uncertainties in the dynamic model, and Monte Carlo analysis is performed. The results show that the method developed in this paper produces tighter trajectory envelopes and smaller terminal state errors without significantly increasing the computational burden when compared with a method that does not penalize sensitivities.
\end{abstract}


\end{frontmatter}

\clearpage
\section{Introduction}
The goal of an optimal control problem (OCP) is to determine the state and control of a dynamical system that optimize a specified performance index subject to dynamic constraints, boundary conditions, and path constraints. In the case of a deterministic OCP, a nominal model is used to compute an optimal reference trajectory. In practice, however, model uncertainties and disturbances create differences between the actual state and the reference state when the reference control is applied. To mitigate these effects, course corrections, known as \textit{guidance updates}, can be used to recompute the reference control at specified update times. Classical guidance frameworks are often described in terms of the following two components: trajectory design and trajectory tracking. Trajectory design employs numerical methods to determine a reference trajectory that satisfies mission requirements and optimizes a specified performance index \cite{bryson-2018,neustadt-2015}. Trajectory tracking involves the design of a feedback control system that reduces the difference between the true state and the reference state, thereby attempting to keep the actual state in close proximity to the reference state \cite{kwakernaak-sivan-1972,franklin-powell-2002,kuo-1987}. From this perspective, trajectory design is typically treated as an open-loop process, while trajectory tracking is a closed-loop process that employs feedback control. The separation between trajectory design and trajectory tracking, however, does not explicitly account for the sensitivity of the optimal solution to model uncertainties. Even when feedback control is used in trajectory tracking, uncertainties in the model can lead to significant performance degradation, including constraint violations, large terminal errors, and control saturation. 

Another approach to guidance is to explicitly take into account uncertainties in the model. Such an approach incorporates sensitivities, where sensitivities quantify how variations in uncertain parameters influence the state and, therefore, system performance. The effects of parameter uncertainties are particularly pronounced in nonlinear systems over long durations where high performance is required. Consequently, incorporating sensitivity information into trajectory design is key to ensure system robustness to model uncertainties. The need for sensitivity-aware trajectory design motivates the development of methods for \textit{desensitized optimal control} (DOC), where DOC seeks to minimize the sensitivity of the nominal trajectory to model uncertainties, producing trajectories that are resilient to both state disturbances and parameter variations \cite{seywald-kumar-2003,shen-seywald-2010,li-peng-2011,seywald-seywald-2024,jawaharlal-taheri-2024,zames-francis-2003,piprek-hong-2022,pilipovsky-tsiotras-2025,robbiani-sagliano-2025,makkapati-dor-2018}. 

Desensitized optimal control methods can be classified as either open-loop or closed-loop. In \textit{open-loop} DOC methods, the reference control is used for guidance. Initial open-loop formulations of desensitized optimal control modified the OCP to explicitly penalize sensitivities to parametric uncertainties. These early approaches introduced sensitivities as augmented states, allowing for the direct penalization of trajectory variations within the objective functional. This formulation enabled control solutions with reduced sensitivity to system disturbances \cite{seywald-kumar-2003}. Subsequent work extended this formulation to constrained OCPs and systems with more complex dynamic models \cite{shen-seywald-2010,li-peng-2011}. To further mitigate the effects of uncertainty, several studies have incorporated feedback into DOC formulations. These \textit{closed-loop} DOC methods introduce feedback gains as additional controls to the OCP such that these gains are optimized along with the reference optimal control \cite{zames-francis-2003,piprek-hong-2022,pilipovsky-tsiotras-2025,robbiani-sagliano-2025}. With the continued application of DOC to systems of increasing complexity, the introduction of sensitivities as state variables has increased the computational burden of solving the resulting OCPs, motivating recent efforts to focus on the development of reduced DOC formulations \cite{seywald-seywald-2024,jawaharlal-taheri-2024,makkapati-dor-2018}. 

Due to the complexity of desensitized optimal control problems, an analytical solution does not often exist and, therefore, numerical methods must be used. Numerical methods for solving OCPs are categorized as either indirect or direct. Indirect methods follow the approach of \textit{optimize then discretize}, where the necessary first-order optimality conditions are derived from the calculus of variations to transcribe the original OCP into a Hamiltonian boundary-value problem (HBVP) \cite{betts-2010,kirk-2004}. Conversely, direct methods follow the approach of \textit{discretize then optimize}. In a direct method, the control and/or state are parameterized, transcribing the continuous-time OCP into a finite-dimensional nonlinear programming problem (NLP) that can be solved by well-developed software \cite{biegler-zavala-2009,gill-murray-2005,byrd-nocedal-2006}. One class of direct methods, known as \textit{direct collocation methods}, is extensively used to compute numerical solutions of continuous-time OCPs \cite{kraft-1985,hargraves-paris-1987,vonstryk-bulirsch-1992,enright-conway-1992,elnagar-kazemi-1995,benson-huntington-2006}. A direct collocation method is a state and control parameterization method in which the dynamics are enforced at a set of specified points called \textit{nodes} or \textit{collocation points}.  Among these methods, \textit{Gaussian quadrature collocation} has received significant attention due to its high accuracy. The most well-developed Gaussian quadrature methods employ Legendre-Gauss (LG) \cite{benson-huntington-2006,rao-benson-2010}, Legendre-Gauss-Radau (LGR) \cite{kameswaran-biegler-2008,garg-patterson-2010,garg-hager-2011,garg-patterson-2011}, or Legendre-Gauss-Lobatto (LGL) points \cite{elnagar-kazemi-1995} where the collocation points are the corresponding roots of Legendre polynomials or linear combinations thereof. When the solution of the optimal control problem is smooth, LG and LGR methods exhibit exponential convergence \cite{hager-hou-2016,hager-hou-2017,hager-liu-2018,chen-du-2019,hager-hou-2019}.

The development of computationally efficient and accurate numerical methods for solving OCPs made it possible to develop \textit{computational optimal guidance and control} methods \cite{lu-2017,dueri-acikmese-2017,scharf-acikmese-2017,lu-brunner-2017,ferranti-keviczky-2017}. These methods recompute the reference optimal control online using updated state information, leveraging onboard optimization to attenuate trajectory disturbances. The development of entry guidance algorithms for onboard use relies on fast and accurate computation of solutions to OCPs. However, existing computational optimal guidance and control approaches solve \textit{nominal} OCPs and do not explicitly account for sensitivities to model uncertainties when recomputing the reference optimal control. As a result, the recomputed optimal trajectories remain sensitive to parametric uncertainties and limit the robustness of the guidance algorithm even when guidance updates are performed frequently. A key goal of this research is to develop a computational method for desensitized optimal guidance that embeds desensitized optimal control into a shrinking-horizon formulation, combining the robustness properties of DOC with computational optimal guidance. 


The contributions of this paper are as follows. First, a computational desensitized optimal guidance method is developed for systems subject to parametric uncertainty. Different from the method of Ref.~\cite{makkapati-dor-2018}, where only a desensitized reference trajectory is obtained, in this paper a desensitized optimal control problem is solved at specified guidance cycle times. Furthermore, at the start of each guidance cycle, the desensitized OCP is the one that corresponds to the remaining horizon. Second, in order to improve computational efficiency, a modified mesh remapping strategy is developed that is different from the method of Ref.~\cite{dennis-hager-2019}. Specifically, this modified mesh remapping strategy not only truncates and remaps the initial guess of the sensitivity matrix on the remaining horizon, but also shifts the sensitivity matrix such that the initial value of the sensitivity matrix on the remaining horizon is zero.  Together, these elements produce near-optimal guidance commands that are less sensitive to parametric uncertainties, thereby improving robustness without significantly increasing computational complexity. Finally, this paper presents an in-depth development and analysis of the method that provides a comprehensive description of the approach, studies representative problems, and assesses computational efficiency and robustness via a rigorous Monte Carlo analysis.

The remainder of this paper is organized as follows. Section~\ref{section: bolza-OCP} presents the general formulation for the continuous-time Bolza optimal control problem. Section~\ref{section: LGR-collocation} describes the discretization of the Bolza OCP using the $hp-$form of the Legendre-Gauss-Radau orthogonal collocation method, followed by the equivalent nonlinear programming problem. Section~\ref{section: DOG-method-overview} details the desensitized optimal guidance methodology that is the subject of this paper. This section includes the introduction of a general open-loop desensitized optimal control problem formulation followed by the simulation architecture and guidance update procedure to demonstrate how these components are integrated within a unified guidance framework. Section~\ref{section: numerical-examples} demonstrates this novel approach on two numerical examples to highlight the advantages of implementing this method. Finally, Section~\ref{section: conclusion} summarizes the key findings of this work.

\section{Bolza Optimal Control Problem}\label{section: bolza-OCP}
Without loss of generality, consider the following optimal control problem expressed in Bolza form in terms of the elapsed time, $t$, where $t\in[t_0,t_f]$. Determine the state, $\mathbf{x}(t)\in \mathbb{R}^{n_x}$, control, $\mathbf{u}(t)\in \mathbb{R}^{n_u}$, the initial time, $t_0$, and final time, $t_f$, that minimize the objective functional 
\begin{equation}
    \mathcal{J}=\mathcal{M}(\mb{x}(t_0),t_0,\mb{x}(t_f),t_f)+\int_{t_0}^{t_f}\mathcal{L}(\mb{x}(t),\mb{u}(t),t) \dt,
\end{equation}
subject to the dynamic constraints
\begin{equation}
    \dfrac{\text{d}\mb{x}}{\dt}=\mb{f}(\mb{x}(t),\mb{u}(t),t),
\end{equation}
the boundary conditions
\begin{equation}
    \mb{b}_{\min} \leq \mb{b}(\mb{x}(t_0),t_0,\mb{x}(t_f),t_f) \leq \mb{b}_{\max},
\end{equation}
and the inequality path constraints
\begin{equation}
    \mb{c}_{\min} \leq \mb{c}(\mb{x}(t),\mb{u}(t),t) \leq \mb{c}_{\max},
\end{equation}
where $n_x$ is the number of state variables and $n_u$ is the number of control variables. The functions $\mathcal{M},\mathcal{L},\mb{f},\mb{b},\mb{c}$ are defined by the following mappings:
\begin{equation}
\begin{array}{lclcl}
    \mathcal{M} & : & \mathbb{R}^{n_x} \times \mathbb{R} \times \mathbb{R}^{n_x} \times \mathbb{R} & \rightarrow & \mathbb{R},\\
    \mathcal{L} & : & \mathbb{R}^{n_x} \times \mathbb{R}^{n_u} \times \mathbb{R} & \rightarrow & \mathbb{R},\\
    \mb{f} & : & \mathbb{R}^{n_x} \times \mathbb{R}^{n_u} \times \mathbb{R} & \rightarrow & \mathbb{R}^{n_x},\\
    \mb{b} & : & \mathbb{R}^{n_x} \times \mathbb{R} \times \mathbb{R}^{n_x} \times \mathbb{R} & \rightarrow & \mathbb{R}^{n_b},\\
    \mb{c} & : & \mathbb{R}^{n_x} \times \mathbb{R}^{n_u} \times \mathbb{R} & \rightarrow & \mathbb{R}^{n_c},
\end{array}
\end{equation}
where $n_b$ and $n_c$ are the number of boundary conditions and inequality path constraints, respectively. Generally, for numerical stability and convenience, it is preferred to write the optimal control problem in terms of the variable $\tau \in [-1,+1]$, formally known as the mesh time. This is achieved by the affine transformation
\begin{equation}
    t(\tau,t_0,t_f)=\dfrac{t_f-t_0}{2}\tau + \dfrac{t_f+t_0}{2}.
\end{equation}
The transformed continuous-time Bolza optimal control problem is then written as follows. Determine the state, $\mb{x}(\tau)\in\mathbb{R}^{n_x}$, the control, $\mb{u}(\tau)\in \mathbb{R}^{n_u}$, the initial time, $t_0$, and final time, $t_f$, that minimize the objective functional
\begin{equation}\label{eq: transformed-bolza-cost}
    \mathcal{J}=\mathcal{M}(\mb{x}(-1),t_0,\mb{x}(+1),t_f) + \dfrac{t_f-t_0}{2}\int_{-1}^{+1} \mathcal{L}(\mb{x}(\tau),\mb{u}(\tau),t(\tau,t_0,t_f))\text{d} \tau,  
\end{equation}
subject to the constraints
\begin{equation}\label{eq: transformed-bolza-constraints}
\begin{aligned}
    \dfrac{\text{d}\mb{x}}{\text{d} \tau}=&\dfrac{t_f-t_0}{2} \mb{f}(\mb{x}(\tau),\mb{u}(\tau),t(\tau,t_0,t_f)), \\
    \mb{b}_{\min} \leq & \mb{b}(\mb{x}(-1),t_0,\mb{x}(+1),t_f) \leq \mb{b}_{\max}, \\
    \mb{c}_{\min} \leq & \mb{c}(\mb{x}(\tau),\mb{u}(\tau),t(\tau,t_0,t_f)) \leq \mb{c}_{\max}.
\end{aligned}
\end{equation}

\section{Legendre-Gauss-Radau Collocation}\label{section: LGR-collocation}
The continuous-time Bolza optimal control problem is discretized using the $\textit{hp}-$form of the Legendre-Gauss-Radau collocation method as presented in Refs.~\cite{garg-patterson-2010,garg-hager-2011,garg-patterson-2011,darby-hager-2011,darby-hager-rao-2011,patterson-hager-2015,liu-hager-2015,liu-rao-2017}. In the multiple-interval form of the LGR collocation method, the independent variable $\tau\in[-1,+1]$ is partitioned into $K$ mesh intervals such that mesh interval $\mathcal{I}_k$ is defined on $[T_{k-1},T_k]\subseteq[-1,+1]$, $(k = 1,...,K)$. Let $\mb{x}^{(k)}(\tau)$ and $\mb{u}^{(k)}(\tau)$ be the state and control in mesh interval $\mathcal{I}_k$, $(k=1,...,K)$. The state is approximated as
\begin{equation}
    \mb{x}^{(k)}(\tau) \approx \mb{X}^{(k)}(\tau) = \sum_{j=1}^{N_k+1} \mb{X}^{(k)}_j \ell_j^{(k)}(\tau),\quad (k=1,...,K),
\end{equation}
\begin{equation}
    \ell^{(k)}_j(\tau) = \prod_{\substack{l=1\\l\neq j}}^{N_k+1} \dfrac{\tau - \tau_l^{(k)}}{\tau_j^{(k)}-\tau_l^{(k)}},\quad (k=1,...,K),
\end{equation}
where $\ell_j^{(k)}(\tau)$, $(j=1,...,N_k+1)$ is a basis of Lagrange polynomials on mesh interval $\mathcal{I}_k$, ($\tau_1^{(k)},...,\tau_{N_k}^{(k)}$) are the $N_k$ LGR collocation points on the interval $[T_{k-1},T_k)$, and $\tau^{(k)}_{N_k+1}=T_k$ is a noncollocated support point. The derivative of the state approximation is then
\begin{equation}
    \dfrac{\text{d}\mb{x}^{(k)}(\tau)}{\text{d}\tau} \approx
    \dfrac{\text{d}\mb{X}^{(k)}(\tau)}{\text{d}\tau}=\sum_{j=1}^{N_k+1} \mb{X}_j^{(k)}\dfrac{\text{d} \ell_j^{(k)}(\tau)}{\text{d}\tau}.
\end{equation}
Now, let $\mb{X}_j^{(k)}$ denote the state approximation at $\tau_j^{(k)}$, $(j=1,...,N_k+1)$ and $\mb{U}_i^{(k)}$ denote the control approximation at $\tau^{(k)}_i$, $(i=1,...,N_k)$. The multiple-interval LGR collocation defect constraints are determined by
\begin{equation}\label{eq: defect-constraints}
    \sum_{j=1}^{N_k+1}[D_{ij}^{(k)}]^\text{LGR} \mb{X}_j^{(k)} - \dfrac{t_f - t_0}{2} \dfrac{T_k-T_{k-1}}{2}\mb{f}(\mb{X}^{(k)}_i,\mb{U}^{(k)}_i,\tau^{(k)}_i;t_0,t_f)=\mb{0},\ \ \ \ (i=1,...,N_k),
\end{equation}
where
\begin{equation}
    [D_{ij}^{(k)}]^\text{LGR} \equiv \dfrac{\text{d}\ell_j^{(k)}(\tau^{(k)}_i)}{\text{d}\tau},\quad (i=1,...,N_k),\quad (j=1,...,N_k+1).
\end{equation}
denotes the element in the $i^{th}$ row and $j^{th}$ column of the $N_k \times (N_k+1)$ LGR differentiation matrix in mesh interval $\mathcal{I}_k$. Continuity in the state at the interior mesh points is enforced by
\begin{equation}\label{eq: continuity-condition}
    \mb{X}_{N_k+1}^{(k)} = \mb{X}_1^{(k+1)},\quad (k=1,...,K-1).
\end{equation}
Lastly, the total number of collocation points, $N$, is determined by the summation $N = \sum_{k=1}^{K} N_k$. Approximating the integral in Eq.~\eqref{eq: transformed-bolza-cost} via Gauss quadrature and the constraints in Eq.~\eqref{eq: transformed-bolza-constraints} at the $N_k$ LGR nodes leads to the following nonlinear programming problem. Minimize the cost function
\begin{equation}
\begin{aligned}
    \mathcal{J} \approx& \mathcal{M}(\mb{X}_1^{(1)},\tau_1^{(1)},\mb{X}^{(K)}_{N_k+1},\tau_{N_k+1}^{(K)};t_0,t_f)\\
    &+\dfrac{t_f-t_0}{2}\sum_{k=1}^K \dfrac{T_k - T_{k-1}}{2}
    \sum_{i=1}^{N_k} w_i^{(k)}\mathcal{L}(\mb{X}_i^{(k)},\mb{U}_i^{(k)},\tau_i^{(k)};t_0,t_f),
\end{aligned}
\end{equation}
subject to Eq.~\eqref{eq: defect-constraints}, the continuity condition in Eq.~\eqref{eq: continuity-condition}, and the constraints
\begin{equation}\label{eq: NLP-constraints}
\begin{aligned}
      \mb{b}_\text{min} \leq & \mb{b}(\mb{X}_1^{(1)},\tau_1^{(1)},\mb{X}_{N_K+1}^{(K)},\tau^{(K)}_{N_K+1}) \leq \mb{b}_\text{max}, \\
     \mb{c}_\text{min} \leq & \mb{c}(\mb{X}_i^{(k)},\mb{U}_i^{(k)},\tau_i^{(k)};t_0,t_f) \leq \mb{c}_\text{max},
\end{aligned}
\end{equation}
$\forall (i=1,...,N_k),(k=1,...,K)$ where $w_i^{(k)}$, $(i=1,...,N_k)$ are the LGR quadrature weights in mesh interval $\mathcal{I}_k$. The resulting NLP given in Eq.~\eqref{eq: defect-constraints}-\eqref{eq: NLP-constraints} can be solved using well-known NLP solvers that exploit the natural sparsity of the problem. These solvers include but are not limited to Interior Point OPTimizer (IPOPT) \cite{biegler-zavala-2009}, Sparse Nonlinear OPTimizer (SNOPT) \cite{gill-murray-2005}, and Nonlinear Interior-point Trust Region Optimization (KNITRO) \cite{byrd-nocedal-2006}. 

\section{Method for Desensitized Optimal Guidance}\label{section: DOG-method-overview}
In this section, a computational method for desensitized optimal guidance (DOG) is presented. This method is for application to dynamical systems subject to parameter uncertainties in the dynamic model. The DOG method consists of two cyclic stages. The first stage, described in further detail in Section~\ref{section: DOC-formulation}, determines a reference trajectory that is less sensitive to model uncertainties through the introduction of state sensitivities as augmented state variables. The second stage, presented in Section~\ref{section: Mesh-remapping}, employs guidance updates to recompute the desensitized optimal control at prescribed times along the perturbed trajectory, where computational efficiency is preserved through a mesh remapping strategy that accelerates convergence of the NLP. Specifically, the control solution computed at the previous guidance update is used to initialize the guess for the subsequent update. The cyclic nature of the algorithm is demonstrated through the iterative repetition of stage one and stage two until the trajectory terminates, resulting in a sequence of guidance updates that continually works to mitigate the effects of parametric uncertainties. 

\subsection{Desensitized Optimal Control Formulation}\label{section: DOC-formulation}
This section summarizes the reduced desensitized optimal control formulation developed by Ref.~\cite{makkapati-dor-2018}. This formulation is adapted within the proposed guidance method to generate reference trajectories that are robust to the model imperfections represented by the set of uncertain parameters $\mb{p}\in \mathbb{R}^{n_p}$. The sensitivities of the state with respect to these parameters are the elements of the $n_x \times n_p$ time-varying sensitivity matrix, $\mb{S}(t)$, where $n_x$ and $n_p$ denote the number of state variables and the number of uncertain parameters, respectively. The sensitivity matrix is the solution to the linear, time-varying initial value problem
\begin{equation}\label{eq: sensitivity-IVP}
\begin{array}{lcc}
    \dot{\mb{S}}(t)=\mb{A}(t)\mb{S}(t) + \mb{B}(t), & \mb{S}(t_0) = \mb{0},
\end{array}
\end{equation}
where the matrix $\mb{A}(t)\in \mathbb{R}^{n_x \times n_x}$ is the partial derivative of the continuous dynamics with respect to the state, and the matrix $\mb{B}(t)\in \mathbb{R}^{n_x \times n_p}$ is the partial derivative of the continuous dynamics with respect to the parameters given by
\begin{equation}\label{eq: A(t)-B(t)-definitions}
\begin{aligned}
    \mb{A}(t) & =  \dfrac{\partial \mb{f}}{\partial \mb{x}}(\mb{x}(\mb{p},t),\mb{u}(t),\mb{p},t),\\
    \mb{B}(t) & =  \dfrac{\partial \mb{f}}{\partial \mb{p}}(\mb{x}(\mb{p},t),\mb{u}(t),\mb{p},t).
\end{aligned}
\end{equation}
Note that Eq.~\eqref{eq: A(t)-B(t)-definitions} is evaluated using the nominal value of the parameter. Each element of the sensitivity matrix is introduced as an augmented state variable to the OCP formulation. The decision of which state sensitivities to minimize is dictated by the user-defined output vector $\mb{y}=\mb{g}(\mb{x})$ which is a function of the state and explicitly penalizes trajectory variations within the objective functional. Using the first-order approximation described in Ref.~\cite{makkapati-dor-2018}, the augmented cost
\begin{equation}\label{eq: JA}
    \mathcal{J}_a = \mathcal{J} + \mathbb{E}\left[||\delta \mb{y}(t_f)||^2 + \int_{t_0}^{t_f}||\delta \mb{y}(t)||^2\text{d}t\right],
\end{equation}
is defined in terms of the stochastic error term
\begin{equation}\label{eq: error}
    \mathbb{E}(||\delta \mb{y}(t)||^2) \approx  \text{tr}[ \mb{G} \mb{S}(t) \mb{P} \mb{S}(t)^\top \mb{G}^\top],
\end{equation}
where $\mb{P}$ is a known positive semi-definite parameter covariance matrix, $\mb{G}=\partial \mb{g}/\partial \mb{x}$ is the Jacobian of the output function evaluated using the nominal state, and $\text{tr}[\cdot]$ indicates the matrix trace operator. 
The augmented objective functional is rewritten as a weighted combination of the normalized performance and sensitivity penalties such that
\begin{equation}\label{eq: normalized-cost}
    \bar{\mathcal{J}}_a = w_1\frac{\mathcal{J}}{\mathcal{J}_\text{ref}} + w_2 \frac{\mathcal{J}_f}{\mathcal{J}_{f,\text{ref}}} + w_3 \frac{\mathcal{J}_r}{\mathcal{J}_{r,\text{ref}}},
\end{equation}
where $\mathcal{J}_f$ is the terminal sensitivity penalty, $\mathcal{J}_r$ is the running sensitivity penalty, and the reference values $\mathcal{J}_\text{ref}$, $\mathcal{J}_{f,\text{ref}}$, and $\mathcal{J}_{r,\text{ref}}$ are the associated cost components evaluated on the nominal solution without desensitization. The weights
\begin{equation}\label{eq: sensitivity-weights}
\begin{aligned}
    w_1 =& 1-w_2- w_3,\quad \\
    w_2 =& \beta, \\
    w_3 =& \xi \beta,
\end{aligned}
\end{equation}
are defined in terms of the desensitization parameters $\beta$ and $\xi$, where $\beta$ determines the magnitude of the terminal sensitivity penalty, and $\xi$ specifies the fraction of the weight, $\beta$, assigned to the running sensitivity penalty. Equation~\eqref{eq: sensitivity-weights} satisfies $\beta (1+\xi) \leq 1$ such that $\sum_{i=1}^3 w_i=1$. As a result of this normalization, each cost component is $\mathcal{O}(1)$, allowing the desensitization parameters $\beta$ and $\xi$ to directly reflect a trade-off between performance and robustness. 

The resulting augmented OCP is then given as follows. Determine the state, $\mb{x}(\mb{p},t)\in \mathbb{R}^{n_x}$, control, $\mb{u}(t)\in\mathbb{R}^{n_u}$, and sensitivities, $\mb{S}(t)\in \mathbb{R}^{n_x \times n_p}$, where $t \in[t_0,t_f]$ that minimize the objective functional in Eq.~\eqref{eq: normalized-cost}, subject to the dynamic constraints 
\begin{equation}
\begin{aligned}
\dfrac{\text{d}\mb{x}(\mb{p},t)}{\text{d}t} =& \mb{f}(\mb{x}(\mb{p},t),\mb{u}(t),\mb{p},t),\\
    \dfrac{\text{d}\mb{S}(t)}{\text{d}t} =& \mb{A}(t) \mb{S}(t) + \mb{B}(t),
\end{aligned}
\end{equation}
the boundary conditions
\begin{equation}
   \mb{b}_\text{min}   \leq  \mb{b}(\mb{x}(\mb{p},t_0),t_0,\mb{x}(\mb{p},t_f),t_f) \leq \mb{b}_\text{max}, 
\end{equation}
\begin{equation}
    \mb{S}(t_0)  = \mb{0}, 
\end{equation}
and the inequality path constraints
\begin{equation}
    \mb{c}_\text{min} \leq \mb{c}(\mb{x}(\mb{p},t),\mb{u}(t),\mb{p},t) \leq \mb{c}_\text{max}.
\end{equation}
The solution to this augmented OCP is a reference trajectory that is less sensitive to parametric uncertainties. However, this desensitized OCP formulation treats the control as a fixed function of time such that the control does not change in response to state perturbations. The next section details how the proposed guidance framework introduces control corrections in the form of guidance updates.

\subsection{Shrinking-Horizon Guidance with Mesh Remapping}\label{section: Mesh-remapping}
\begin{figure}[h!]
\begin{center}
	\includegraphics[scale=.82]{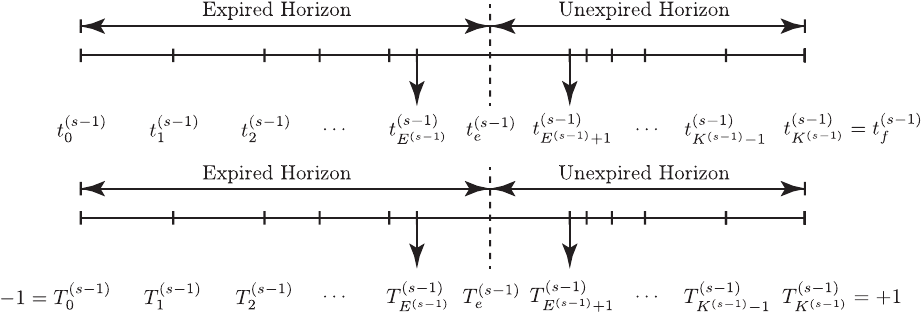}
	\caption{Partitioned mesh for guidance cycle $(s-1)$ \cite{dennis-hager-2019}.}
	\label{fig: mesh-remapping}
\end{center}
\end{figure}

This section presents a shrinking-horizon guidance framework that embeds open-loop desensitized optimal control within successive guidance updates. At each guidance update, the reference solution is recomputed on the remaining trajectory, enabled by a mesh truncation and remapping procedure that leverages LGR node placement of the previous control solution to construct the initial mesh. The novelty of this approach lies in the integration of a reduced desensitized optimal control formulation with a mesh remapping technique to efficiently perform mid-course guidance corrections. Let $s$, $(s=1,...,S)$ denote the current guidance cycle iteration where $\Delta t$ is the guidance cycle duration. The current guidance cycle iteration is then defined for $t \in [t_0^{(s)},t_e^{(s)}]$ where
\begin{equation}
\begin{aligned}
    t_0^{(s)} =& t_0 + s \Delta t,\\
    t_e^{(s)} =& t_0 + (s+1) \Delta t. 
\end{aligned}
\end{equation}
Suppose now that the desensitized reference optimal control is implemented on guidance cycle $s$ where the dynamics subject to parameter uncertainties are given as
\begin{equation}
    \dfrac{\text{d}\tilde{\mb{x}}(\tilde{\mb{p}},t)}{\text{d}t} = \tilde{\mb{f}}(\tilde{\mb{x}}(\tilde{\mb{p}},t),\mb{u}^*(t),\tilde{\mb{p}},t),
\end{equation}
where $\tilde{\mb{f}}\in \mathbb{R}^{n_x}$ represents the true dynamics with model uncertainties, $\tilde{\mb{x}}(\tilde{\mb{p}},t)$ is the simulated state, $\tilde{\mb{p}}$ is the vector of perturbed parameters, and $\mb{u}^*(t)$ is the desensitized reference optimal control. It is important to note that the optimal state, $\mb{x}^*(\mb{p},t)$, does \textit{not} equal the simulated state, $\tilde{\mb{x}}(\tilde{\mb{p}},t)$. The true state deviates from the reference trajectory due to the presence of uncertainties in the dynamic model. At the end of the guidance cycle (denoted by the terminal time $t_e^{(s)}$), the mesh is partitioned into expired and unexpired segments. The expired horizon is then removed, and the initial conditions are updated to the current state. To determine the starting mesh for the next guidance update, the location of the terminal time from the previous guidance cycle, $t_e^{(s-1)}$, must be considered. The terminal time of the previous guidance cycle does not exactly occur at a singular mesh point, but instead occurs between two mesh points as shown in Fig.~\ref{fig: mesh-remapping}. 
To delete the expired mesh, the last mesh point of the expired horizon, $T_{E^{(s-1)}}^{(s-1)}$, must be relocated to the point on the mesh corresponding to the terminal time of the previous guidance cycle, i.e., $T_{E^{(s-1)}}^{(s-1)}=T_e^{(s-1)}$. Based on this notation, the first mesh point of the unexpired horizon is then denoted by $T_{E^{(s-1)}+1}^{(s-1)}$. It is important to emphasize that $T_e^{(s-1)}$ is not a mesh point; therefore, the selection of the guidance cycle duration is independent of the mesh partitioning. The initial conditions are updated as
\begin{equation}\label{eq:reset-at-GC-start}
    \mb{x}_0 = \mb{x}(\mb{p},t_0^{(s)}) = \tilde{\mb{x}}(\tilde{\mb{p}},t_e^{(s-1)}),\quad
    \mb{S}_0 = \mb{S}(t_0^{(s)}) = \mb{0}.
\end{equation}
It is seen in Eq.~\eqref{eq:reset-at-GC-start} that  the sensitivity matrix is re-initialized to zero at the start of a guidance cycle.  The reason for resetting the sensitivity matrix to zero is that the sensitivities are assumed to accumulate across a single guidance cycle as opposed to accumulating across the full horizon.  Furthermore, dependence of the sensitivity matrix is assumed to be taken into account by the state estimator that may be used (in the case where the state is not known perfectly). This approach decouples the uncertainty due to model parameters from any uncertainty in the state.


The newly truncated mesh must then be remapped such that $t \in [t_0^{(s)},t_f]$ maps to $\tau\in[-1,+1]$. This mesh truncation and remapping technique is the result of the work done by Ref.~\cite{dennis-hager-2019}. To remap the mesh corresponding to mesh points $(T_e^{(s-1)},...,T_{K^{(s-1)}}^{(s-1)})$ to $\tau\in[-1,+1]$, the following transformation is used
\begin{equation}
    T_k^{(s)} = \dfrac{2\left(T_{k+E^{(s-1)}-1}^{(s-1)} - T_{E^{(s-1)}}^{(s-1)}\right) - \left(1 - T_{E^{(s-1)}}^{(s-1)} \right)}{1-T_{E^{(s-1)}}^{(s-1)}},\ \ \ \ (k=1,...,K^{(s-1)}-E^{(s-1)}+1).
\end{equation}
The number of LGR collocation points used in each mesh interval on the new mesh is determined by
\begin{equation}
    N_k^{(s)}= N_{k+E^{(s-1)}-1}^{(s-1)},\ \ \ \ (k=1,...,K^{(s-1)}-E^{(s-1)}+1).
\end{equation}
The mesh points on the newly remapped mesh (with the exception of the first mesh point $T_{E^{(s-1)}}^{(s-1)}$) align with those from the previous mesh corresponding to guidance cycle $s-1$. Similarly, the LGR nodes on the mesh to be used in guidance update $s$ are exactly aligned with the LGR nodes from guidance cycle $s-1$ with the exception of the first mesh interval due to the relocation of the first mesh point. The mesh for guidance cycle $s$ then consists of $K^{(s)}=K^{(s-1)}-E^{(s-1)}+1$ mesh intervals where the mesh intervals $\mathcal{I}^{(s)}=[T_{k-1}^{(s)},T_k^{(s)}]$, $(k=1,...,K^{(s)})$ have $N_k^{(s)}$ LGR collocation points. The result of this mesh truncation and remapping procedure is a computationally efficient guidance algorithm that recomputes the optimal trajectory on a reduced mesh with well-placed collocation points for rapidly determining a solution to the nonlinear programming problem.

\subsection{Initial Guess of Sensitivity Trajectory on Remaining Horizon}
\begin{figure}[H]
\begin{center}	\includegraphics[scale=.9]{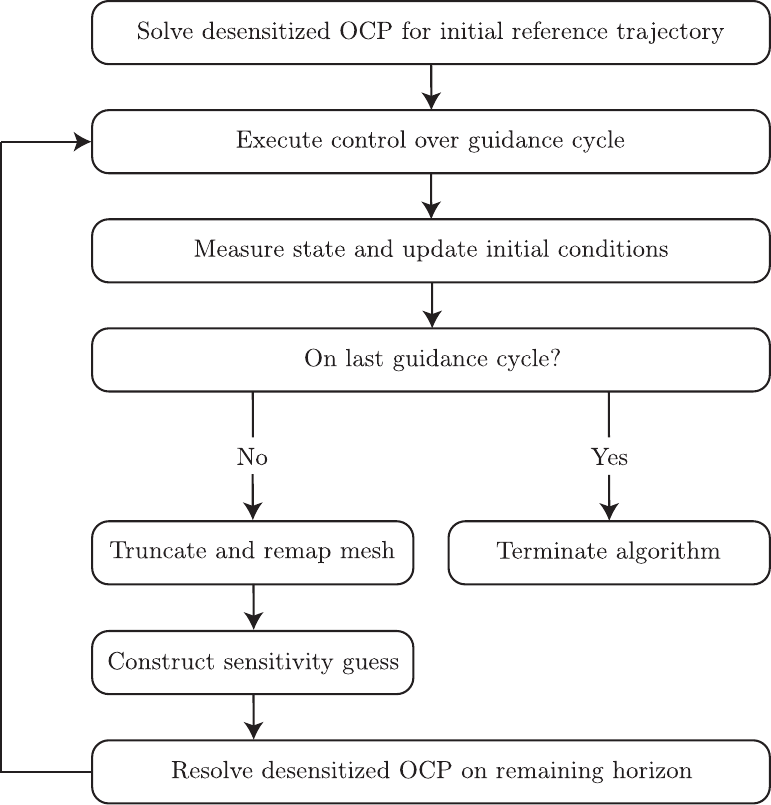}
	\caption{Algorithmic flow diagram for desensitized optimal guidance.}
	\label{fig: dog-flow-chart}
\end{center}
\end{figure}

Different from the method of Ref.~\cite{dennis-hager-2019}, the method in this paper includes sensitivities, and these sensitivities are re-initialized to zero at the start of each guidance cycle (see Section  \ref{section: Mesh-remapping}). Consequently, in addition to providing an initial guess of the state and control on the remaining horizon at the start of a guidance cycle, an initial guess must also be supplied for the sensitivity matrix on the remaining horizon. Moreover, as described in Section \ref{section: Mesh-remapping} and shown in Eq.~\eqref{eq:reset-at-GC-start}, the sensitivity matrix must be re-initialized such that the initial value of the sensitivity matrix is zero at the start of a guidance cycle. This reset is achieved using the following two-step procedure: (1) the mesh truncation and remapping is applied to the sensitivity matrix on the remaining horizon; and (2) the value of the sensitivity matrix at the start of the remaining horizon (after truncation and remapping) is subtracted from the sensitivity matrix obtained in step (1).  Quantitatively, this two-step procedure looks as follows. First, suppose that $\m{S}^{(s)}(t)$ is the truncated and remapped sensitivity matrix for guidance cycle $s$. The initial guess of the sensitivity matrix for guidance cycle $s$ is then constructed 
\begin{equation}
    \hat{\mb{S}}^{(s)}(t) = \mb{S}^{(s)}(t) -\mb{S}^{(s)}(t_0^{(s)}).
\end{equation}
The aforementioned re-initialization described in this section essentially shifts the sensitivity matrix $\mb{S}^{(s)}(t)$ by $\mb{S}^{(s)}(t_0)$, thereby retaining the shape of the previously computed sensitivity  trajectory while making
$\hat{\mb{S}}^{(s)}(t_0)=\m{0}$. The resulting remapped sensitivity trajectory provides an appropriate  initialization $\m{S}(t)$ on the remaining horizon, thereby reducing the computational effort required by the NLP solver. The algorithmic flow of this method can be observed in Fig.~\ref{fig: dog-flow-chart}. 

\section{Numerical Examples}\label{section: numerical-examples}
The following examples demonstrate the efficacy of desensitized optimal guidance in practice. Section~\ref{section: zermelo-problem} solves Zermelo's navigation problem, a commonly studied optimal control problem for desensitized trajectory optimization techniques \cite{seywald-seywald-2024}. Section~\ref{section: rlv-entry} features the unconstrained reusable launch vehicle (RLV) entry problem, offering a more complex real-world application \cite{betts-2010}. Through Monte Carlo analysis, each example illustrates the ability of the framework to mitigate dynamic model uncertainties and, as a result, improve system performance by producing trajectories that achieve realized performance indices closer to the true optimum. The desensitized optimal reference trajectories are computed using the general purpose optimal control software, $\mathbb{GPOPS-II}$ \cite{patterson-rao-2014}. Mesh refinement is permitted with a maximum of 10 iterations and a mesh tolerance of $1 \times 10^{-5}$ \cite{liu-rao-2017}. The NLP solver, IPOPT, is used with an NLP tolerance of $1 \times 10^{-7}$ and the maximum allowable iterations set to 2000. The first and second NLP derivatives were supplied by the open-source algorithmic differentiation software $\textit{ADiGator}$ \cite{weinstein-rao-2017}. The same mesh and NLP tolerance was used to compute both the initial and updated desensitized optimal trajectories. Finally, all simulations were computed using MATLAB R2025a on a 2025 Apple M3 Ultra Mac Studio equipped with 28 cores and 256 GB of RAM. 

\subsection{Zermelo's Navigation Problem}\label{section: zermelo-problem}
\begin{figure}[h!]
\centering
\includegraphics[width=.5\textwidth]{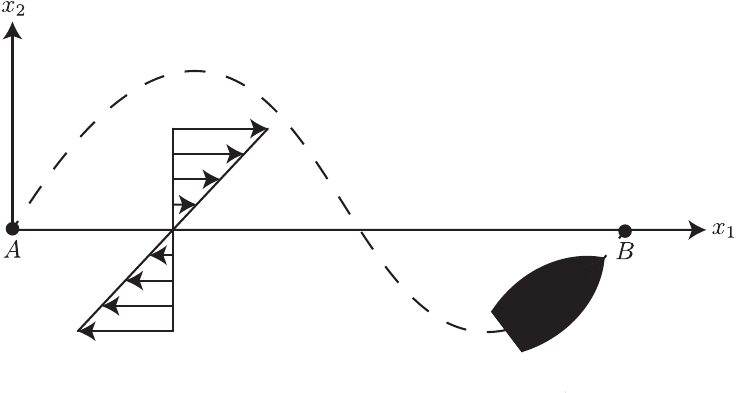}
\caption{Zermelo's navigation problem.}
\label{fig: zermelo-diagram}
\end{figure}

Consider Zermelo's navigation problem which features a river flowing in the $x_1$-direction and a boat traveling downstream. The shore lies parallel to the $x_1$-axis such that the river current increases linearly with the distance from the river center, that is, $v_{current}=c \cdot x_2$ where $c$ is a constant parameter. 
A schematic of the system is provided in Fig.~\ref{fig: zermelo-diagram}. 

\begin{figure}[h!]
    \centering
    \subfloat[Optimal state\label{fig: zermelo-state}]{
        \includegraphics[scale=0.3]{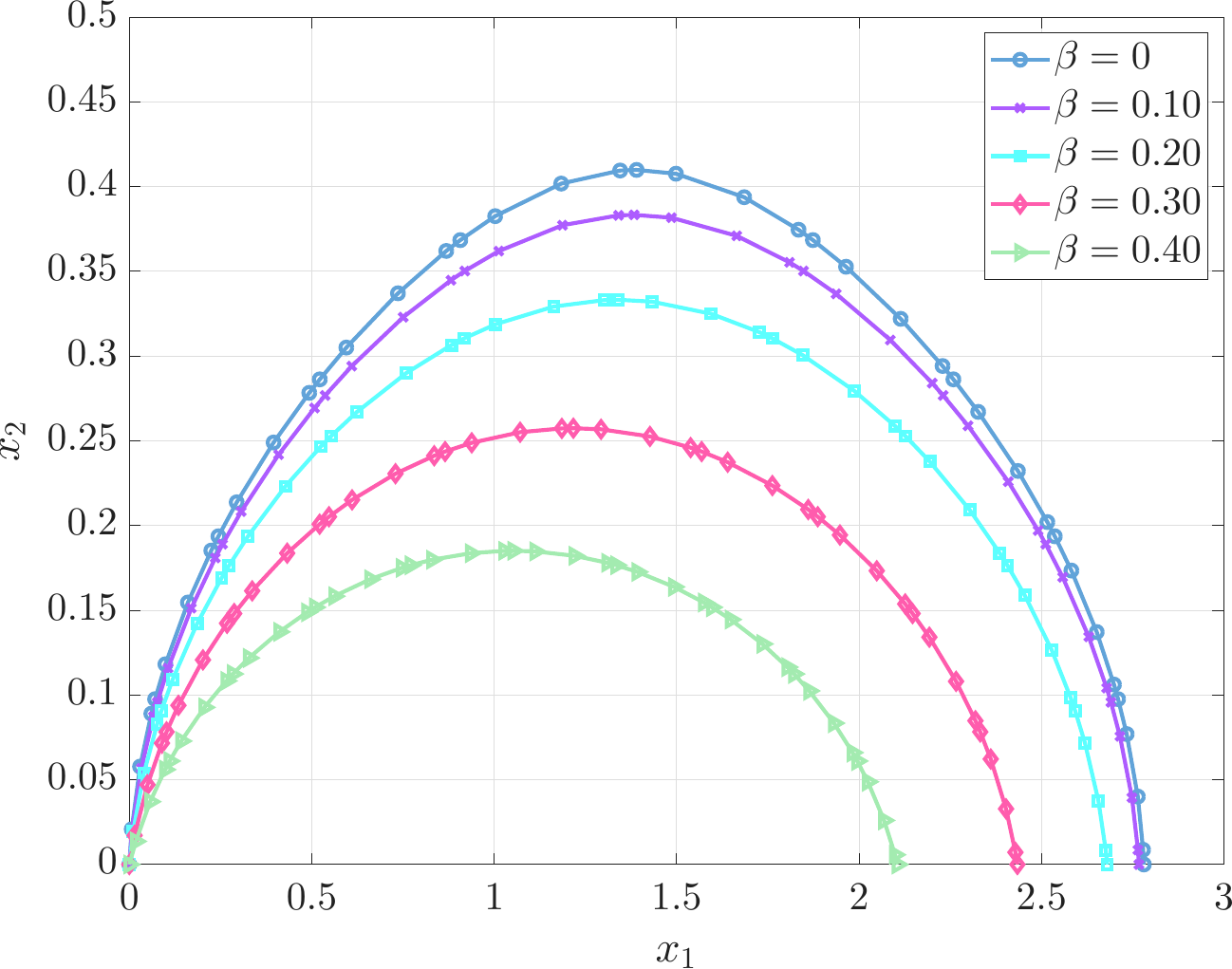}}
    \hspace{0.02\linewidth}
    \subfloat[Optimal control\label{fig: zermelo-control}]{
        \includegraphics[scale=0.3]{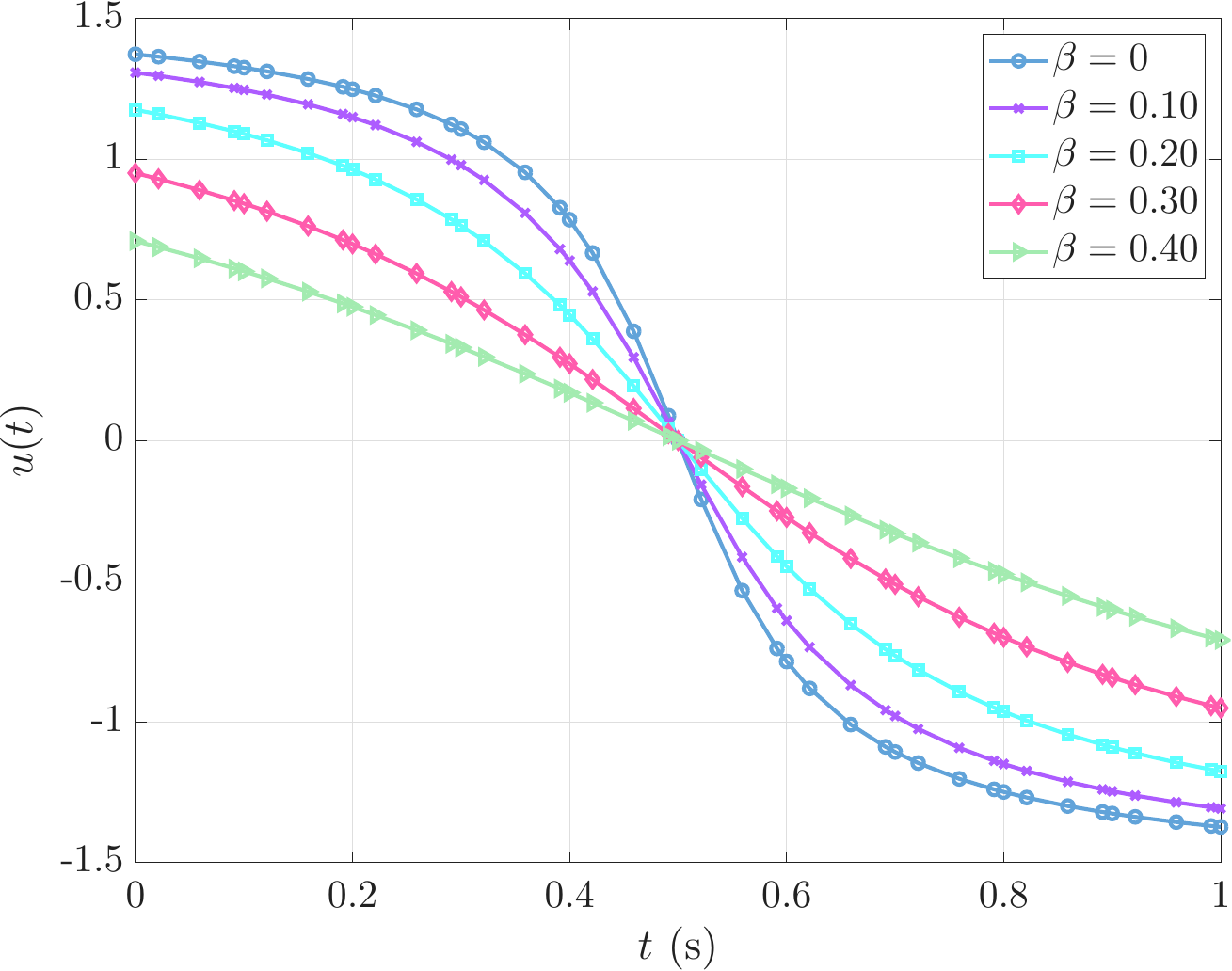}}\\
    \subfloat[Sensitivity of $x_1$ with respect to $c$\label{fig: zermelo-sensitivity}]{
        \includegraphics[scale=0.3]{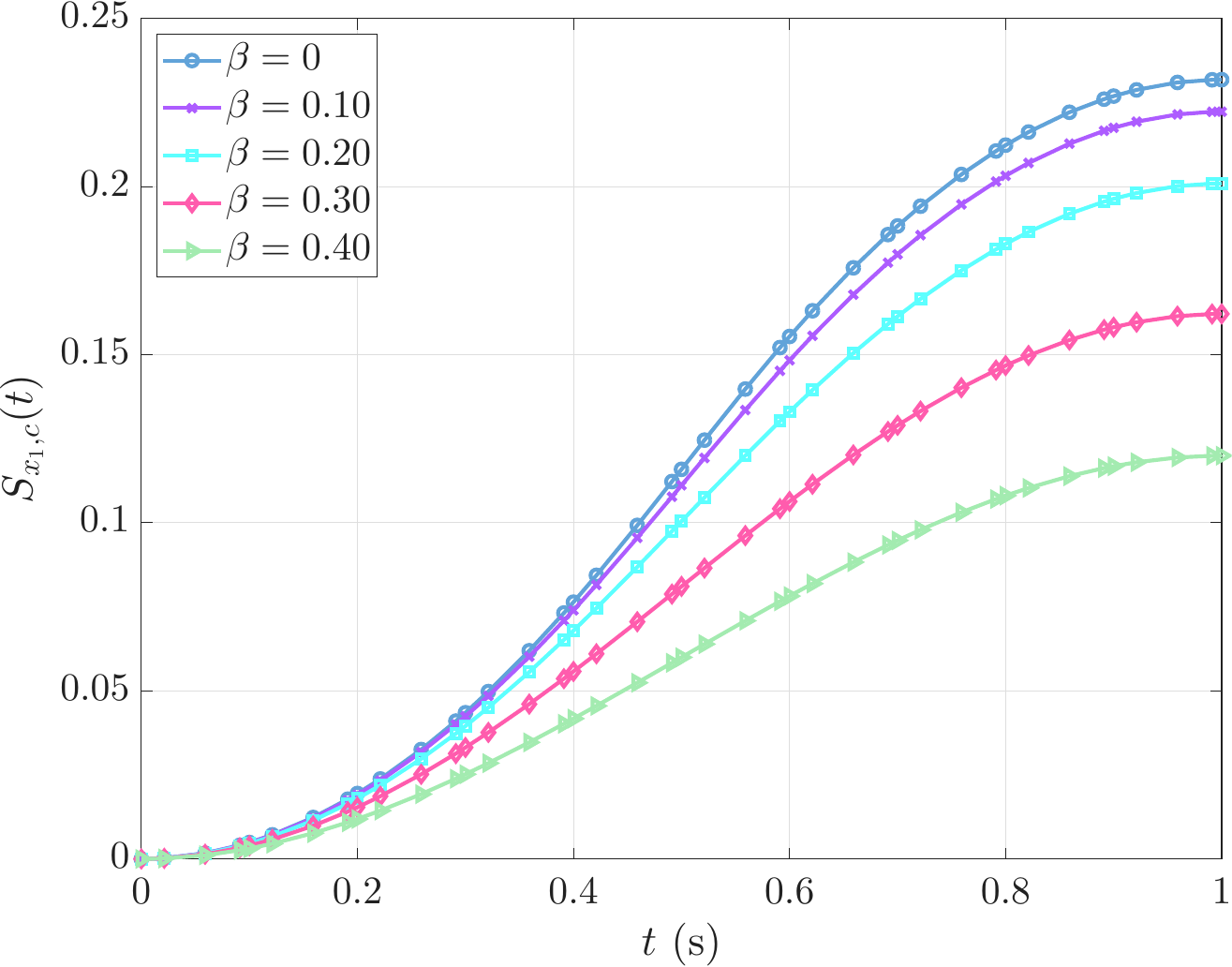}}
    \caption{The reference solution for Zermelo's navigation problem with $\beta\in \{0,0.10,0.20,0.30,0.40\}$.}
    \label{fig: zermelo-reference}
\end{figure}

\noindent Determine the state, $\mb{x}(t)$, and control, $\mb{u}(t)$, that minimize the objective functional
\begin{equation}\label{eq: zermelo-cost}
    \mathcal{J} = -x_1(t_f),
\end{equation}
subject to the dynamic constraints
\begin{equation}\label{eq: zermelo-dynamics}
\begin{array}{rcl}
    \dot{x}_1 & = & \cos(u) + c x_2,\\
    \dot{x}_2 & = & \sin(u), 
\end{array}
\end{equation}
and the boundary conditions
\begin{equation}\label{eq: zermelo-boundary-conditions}
\begin{array}{rclcl}
    x_1(t_0) & = & x_1(0) & = & 0,\\
    x_2(t_0) & = & x_2(0) & = & 0,\\
    x_2(t_f) & = & x_2(1) & = & 0.\\
\end{array}
\end{equation}
The nominal parameter value was chosen as $c=10$. The goal of this problem is to maximize the horizontal distance traveled by the boat in one second. Suppose, however, that the boat is now subject to model uncertainties, specifically in the parameter $c$. To minimize control sensitivities to these uncertainties, sensitivities of the state with respect to the parameter must be introduced as additional state variables to the OCP formulation. For this particular problem, it is desirable to minimize the distribution of the state, $x_1$, at the final time. This requires including an additional term in the objective functional such that variations in $x_1(t_f)$ are penalized. The normalized augmented objective functional is
\begin{equation}\label{eq: zermelo-augmented-cost}
    \bar{\mathcal{J}}_a = w_1 \frac{\mathcal{J}}{\mathcal{J}_\text{ref}} + w_2 \frac{\mathcal{J}_f}{\mathcal{J}_{f,\text{ref}}},
\end{equation}
where $\mathcal{J}_f = \mathbb{E}(|| \delta y(t_f)||^2)$ and $y(t)=x_1(t)$ is the penalized output. The Jacobian of the mapping $\mb{y}=\mb{g}(\mb{x}(t))$ is then  
\begin{equation}\label{eq: zermelo-jacobian}
    \mathbf{G} = 
    \begin{bmatrix}
        1 & 0
    \end{bmatrix},
\end{equation}
the sensitivity matrix is 
\begin{equation}\label{eq: sensitivity-matrix}
    \mathbf{S} = 
    \begin{bmatrix}
        S_{x_1,c}\\
        S_{x_2,c}
    \end{bmatrix},
\end{equation}
and the parameter covariance $P$ is set such that one standard deviation corresponds to $10$ percent of the nominal parameter value. The sampled parameters maintain a normal distribution with mean value $c$ and covariance $P$, expressed as
\begin{equation}
    \tilde{c}=\mathcal{N}(c,0.1^2c^2),
\end{equation}
where $\tilde{c}$ is representative of a perturbed parameter value. Control corrections were introduced via guidance updates at a rate of 5 Hz. A set of five Monte Carlo simulations was performed, each corresponding to a different terminal desensitization weight, $\beta$, with 1000 trials executed per simulation. The desensitization parameter, $\xi$, was set to zero as for this example, there is no running sensitivity penalty. Note that for $\beta=0$, control desensitization is absent, and the method relaxes to a standard optimal guidance method \cite{dennis-hager-2019}. The computed initial reference state and control are included in Fig.~\ref{fig: zermelo-state} and \ref{fig: zermelo-control}, respectively, for each value of the desensitization weight. Figure~\ref{fig: zermelo-sensitivity} shows the sensitivity of the state, $x_1$, with respect to the uncertain parameter, $c$. 

\begin{figure}[h!]
    \centering
    \subfloat[Probability density function for $x_1(t_f)$\label{fig: zermelo-pdf}]{
        \includegraphics[scale=0.3]{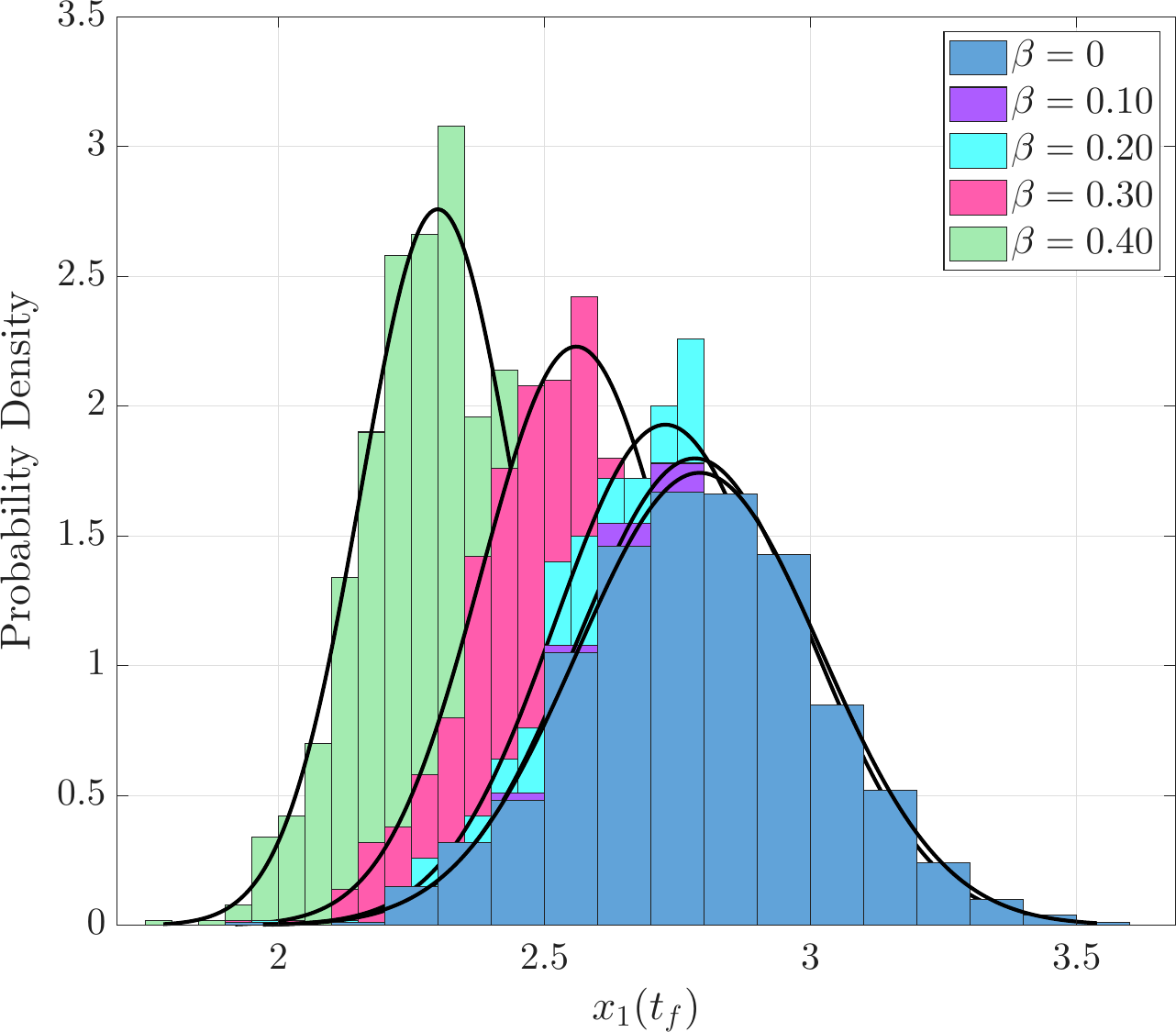}}
    \hspace{0.02\linewidth}
    \subfloat[State trajectory envelope\label{fig: zermelo-envelope}]{
        \includegraphics[scale=0.3]{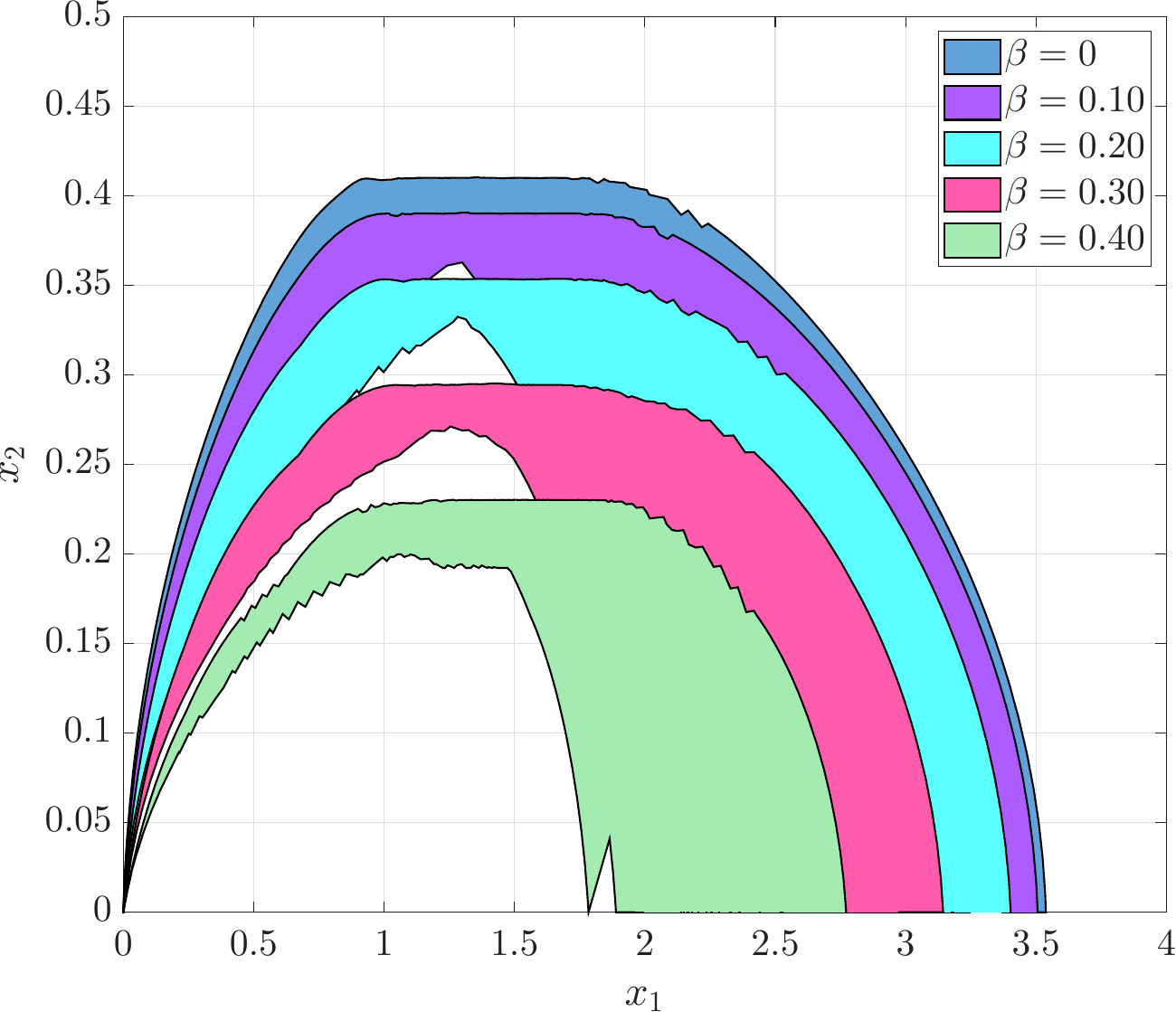}}
    \caption{Results computed from 1000 Monte Carlo trials for each of five values of $\beta\in \{0,0.10,0.20,0.30,0.40\}$.}
    \label{fig: zermelo-MC}
\end{figure}

The results of the Monte Carlo analysis are presented in Fig.~\ref{fig: zermelo-MC}. Figure~\ref{fig: zermelo-pdf} depicts the probability density function of the simulated performance index while Fig.~\ref{fig: zermelo-envelope} displays the state trajectory envelopes. For this particular problem, increasing the desensitization weight produces a more consistent realized performance index, with a reduction in variability achieved at the expense of a modest decrease in the mean optimal cost. Table~\ref{tab: zermelo-trade-study} details the mean achieved terminal horizontal distance, $\mu_{\tilde{x}_{1_f}}$, and corresponding standard deviation, $\sigma_{\tilde{x}_{1_f}}$; the average computational time required for each guidance update denoted $\mu_{CPU}$, and the average percentage of successfully converged guidance updates, $\bar{R}_c$, measured by
\begin{equation}
    \bar{R}_c = \frac{1}{N_{MC}} \sum_{i=1}^{N_{MC}} \frac{n_i^{\text{conv}}}{n_i^{\text{tot}}} * 100,
\end{equation}
where $N_{MC}$ is the number of Monte Carlo trials, $n_i^{\text{conv}}$ is the number of successfully converged guidance updates, and $n_i^{\text{tot}}$ is the total number of attempted guidance updates for the $i^{\text{th}}$ Monte Carlo trial where $i = (1,...,N_{MC})$.
\begin{table}[h!]
\caption{\label{tab: zermelo-trade-study} Zermelo trade study from 1000 Monte Carlo trials for $\beta \in \{0,0.10,0.20,0.30,0.40\}$.}
\centering
\begin{tabular}{ccccccc}
\hline
$\beta$ & $\mu_{\tilde{x}_{1_f}}$ & $\sigma_{\tilde{x}_{1_f}}$ & $\mu_{CPU}$ (s) & $\bar{R}_c$ (\%) & performance loss (\%) & robustness gain (\%)\\
\hline
0.00 & 2.7914 & 0.2288 & 0.2303 & 96.30 & 0 & 0  \\
0.10 & 2.7825 & 0.2218 & 0.0279 & 98.85 & 0.32 & 3.06\\
0.20 & 2.7271 & 0.2068 & 0.2090 & 99.18 & 2.30 & 9.62\\
0.30 & 2.5596 & 0.1789 & 0.1191 & 99.68 & 8.30 & 21.81\\
0.40 & 2.2996 & 0.1446 & 0.1242 & 99.15 & 17.62 & 36.80\\
\hline
\end{tabular}
\end{table}

\clearpage
\subsection{Reusable Launch Vehicle Entry}\label{section: rlv-entry}
Consider the following unconstrained reusable launch vehicle atmospheric entry optimal control problem, augmented from Ref.~\cite{betts-2010}, featuring a vehicle modeled as a point mass in motion over a spherical non-rotating Earth. Determine the state, $\mb{x}(t)$, control, $\mb{u}(t)$, and final time, $t_f$, that maximize the crossrange, mathematically described by the minimization of the objective functional
\begin{equation}\label{eq: rlv-cost}
    \mathcal{J} = -\phi(t_f),
\end{equation}
while subject to the dynamic constraints
\begin{equation}\label{eq: rlv-dynamics}
\begin{split}
\dot{r}       = & v \sin{\gamma},\\
\dot{\theta}  = & \dfrac{v\cos{\gamma}\sin{\psi}}{r \cos{\phi}},\\
\dot{\phi}    = & \dfrac{v \cos{\gamma}\cos{\psi}}{r},\\
\dot{v}       = & -D - g\sin{\gamma} ,\\
\dot{\gamma}  = & \dfrac{L \cos{\sigma}}{v} - \left(\dfrac{g}{v} - \dfrac{v}
                  {r}\right) \cos{\gamma},\\
\dot{\psi}    = & \dfrac{L \sin{\sigma}}{v\cos{\gamma}} - \dfrac{v}{r} 
                  \cos{\gamma}\sin{\psi}\tan{\phi},
\end{split}
\end{equation}
and the boundary conditions
\begin{equation}
\begin{array}{lcrclcr}
    r(t_0) & = & r_0 & , & r(t_f) & = & r_f, \\
    \theta(t_0) & = & \theta_0 & , & \theta(t_f) & = & \theta_f, \\
    \phi(t_0) & = & \phi_0 & , & \phi(t_f) & = & \text{free}, \\
    v(t_0) & = & v_0 & , & v(t_f) & = & v_f, \\
    \gamma(t_0) & = & \gamma_0 & , & \gamma(t_f) & = & \gamma_f, \\
    \psi(t_0) & = & \psi_0 & , & \psi(t_f) & = & \text{free},
\end{array}
\end{equation}
where $r$ is the radial distance from the center of the Earth, $\theta$ is the longitude, $\phi$ is the latitude, $v$ is the speed of the vehicle, $\gamma$ is the flight-path angle, $\psi$ is the azimuth angle, $\alpha$ is the angle of attack, and $\sigma$ is the bank angle. The angle of attack rate and bank angle rate are the control variables, $u_\alpha$ and $u_\sigma$, respectively. The gravitational acceleration is given by $g=\mu/r^2$ where $\mu$ is the gravitational parameter. The radius of the Earth is denoted $R_e$. The lift, $L$, and drag, $D$, are
\begin{equation}\label{eq: rlv-aerodynamic-model}
\begin{aligned}
    L =& \dfrac{qSC_L}{m},\\
    D =& \dfrac{qSC_D}{m},
\end{aligned}
\end{equation}
where the dynamic pressure, $q$, is defined by $q=\rho v^2/2$, $S$ is the surface area of the vehicle, $C_L$ is the lift coefficient, $C_D$ is the drag coefficient, and $m$ is the vehicle mass. The atmospheric density, $\rho$, is modeled as $\rho = \rho_0 \exp (-h/H)$ where $\rho_0$ is the atmospheric density at sea level, $h=r-R_e$ is the altitude, and $H$ is the density scale-height. The aerodynamic model
\begin{equation}\label{eq: rlv-lift-coeff}
    C_L = C_{L_0} + C_{L_1} \alpha,
\end{equation}
\begin{equation}\label{eq: rlv-drag-coeff}
    C_D = C_{D_0} + C_{D_1} \alpha + C_{D_2} \alpha^2,
\end{equation}
is a function of the zero-lift drag coefficient, $C_{D_0}$, through which perturbations will be introduced for this numerical demonstration (see Table~\ref{tab: rlv-relevant-constants} for important model parameter values). Table~\ref{tab: rlv-boundary-conditions} presents the initial and terminal state values.
\begin{table}[b!]
\caption{\label{tab: rlv-relevant-constants} RLV entry relevant constants}
\centering
\begin{tabular}{lccc}
\hline
Description & Symbol & Value & Unit \\
\hline 
Reference area & $S$ & 249.9092 & m$^2$\\
Mass & $m$ & 92079.2526 & kg\\
Equatorial radius of the Earth & $R_e$ & 6.378145$\times 10^6$ & m\\
Sea level atmospheric density & $\rho_0$ & 1.225 & kg/m$^3$\\
Earth gravitational parameter & $\mu$ & 3.986$\times 10^{14}$ & m$^3$/s$^2$\\
Density scale-height & $H$ & 7254.24 & m\\
Zero-lift drag coefficient & $C_{D_0}$ & 0.0785 & -\\
Drag coefficient & $C_{D_1}$ & -0.3529 & - \\
Drag coefficient & $C_{D_2}$ & -0.3529 & - \\
Zero-lift coefficient & $C_{L_0}$ & -0.2070 & -\\
Lift coefficient & $C_{L_1}$ & 1.6756 & -\\
\hline
\end{tabular}
\end{table}
\begin{table}[b!]
   \caption{\label{tab: rlv-boundary-conditions} RLV entry initial and terminal state values}
    \centering 
   \begin{tabular}{lcccc} 
      \hline 
      Description & Symbol & Initial & Terminal & Unit\\
      \hline 
      Time & $t$ & 0 & Free & s\\
      Altitude & $h$ & 79.248 & 24.384 & km\\
      Longitude & $\theta$ & 0 & 75 & deg\\
      Latitude & $\phi$ & 0 & Free & deg\\
      Speed & $v$ & 7.803 & 0.762 & km/s\\
      Flight path angle & $\gamma$ & -1 & -5 & deg\\
      Azimuth angle & $\psi$ & 90 & Free & deg\\
      \hline
   \end{tabular}
\end{table}
For this particular problem, it is desirable to minimize the variation in the longitude, latitude, speed, and flight path angle at the final time in response to the parametric uncertainties introduced through the density scale-height and zero-lift drag coefficient. The penalty function is defined as $\mb{y}=(\theta,\phi,v,\gamma)$. The Jacobian of the mapping $\mb{y}=\mb{g}(\mb{x}(t))$ is then
\begin{equation}
    \mb{G} = 
    \begin{bmatrix}
        0 & 1 & 0 & 0 & 0 & 0\\
        0 & 0 & 1 & 0 & 0 & 0\\
        0 & 0 & 0 & 1 & 0 & 0\\
        0 & 0 & 0 & 0 & 1 & 0\\
    \end{bmatrix},
\end{equation}
and the sensitivity matrix  takes the form
\begin{equation}\label{eq: rlv-sensitivity-matrix}
    \mb{S}(t) = 
    \begin{bmatrix}
        S_{r,H}(t) & S_{r,{C_D}_0}(t)\\
        S_{\theta,H}(t) & S_{\theta,{C_D}_0}(t)\\
        S_{\phi,H}(t) & S_{\phi,{C_D}_0}(t)\\
        S_{v,H}(t) & S_{v,{C_D}_0}(t)\\
        S_{\gamma,H}(t) & S_{\gamma,{C_D}_0}(t)\\
        S_{\psi,H}(t) & S_{\psi,{C_D}_0}(t)\\
    \end{bmatrix},
\end{equation}
with the zero initial condition $\mb{S}(t_0)=\mb{0}$. To determine the nominal solution to this augmented OCP, the angle of attack rate and bank angle rate are implemented as controls to enforce continuity in the angle of attack and bank angle at guidance updates. Each element of the sensitivity matrix in Eq.~\eqref{eq: rlv-sensitivity-matrix} is introduced as a state variable to the OCP. The dynamic constraints are then appended to now include Eq.~\eqref{eq: rlv-dynamics}, the sensitivity dynamics, and the control model
\begin{equation}
\begin{aligned}
    \dot{\alpha} &= u_\alpha,\\
    \dot{\sigma} &= u_\sigma.\\
\end{aligned}
\end{equation}
The parameter covariance, $\mb{P}$, is set such that one standard deviation equates to one percent of the nominal value of the parameter. The perturbed parameter values, $\tilde{\mb{p}}$, are sampled from the Gaussian distribution $\tilde{\mb{p}} \sim \mathcal{N}(\mb{p},\mb{P})$, that is
\begin{equation}
    \begin{bmatrix}
        \tilde{H} \\ \tilde{C}_{D_0}
    \end{bmatrix} \sim
    \mathcal{N} \left( 
    \begin{bmatrix}
        H \\ C_{D_0}
    \end{bmatrix},
    0.01^2
    \begin{bmatrix}
        H^2 & 0 \\
        0 & C^2_{D_0}
    \end{bmatrix}
    \right),
\end{equation}
where $\mb{p}$ denotes the nominal parameter value. During simulation, guidance updates are applied at 20-second intervals. Note the angle of attack and bank angle reference solutions are indirectly controlled through implementation of their respective rates as the controls for this optimal control problem formulation. In all guidance simulations, the computed angle of attack and bank angle are then used as the controls. 

\begin{figure}[b!]
\centering
\includegraphics[scale=0.3]{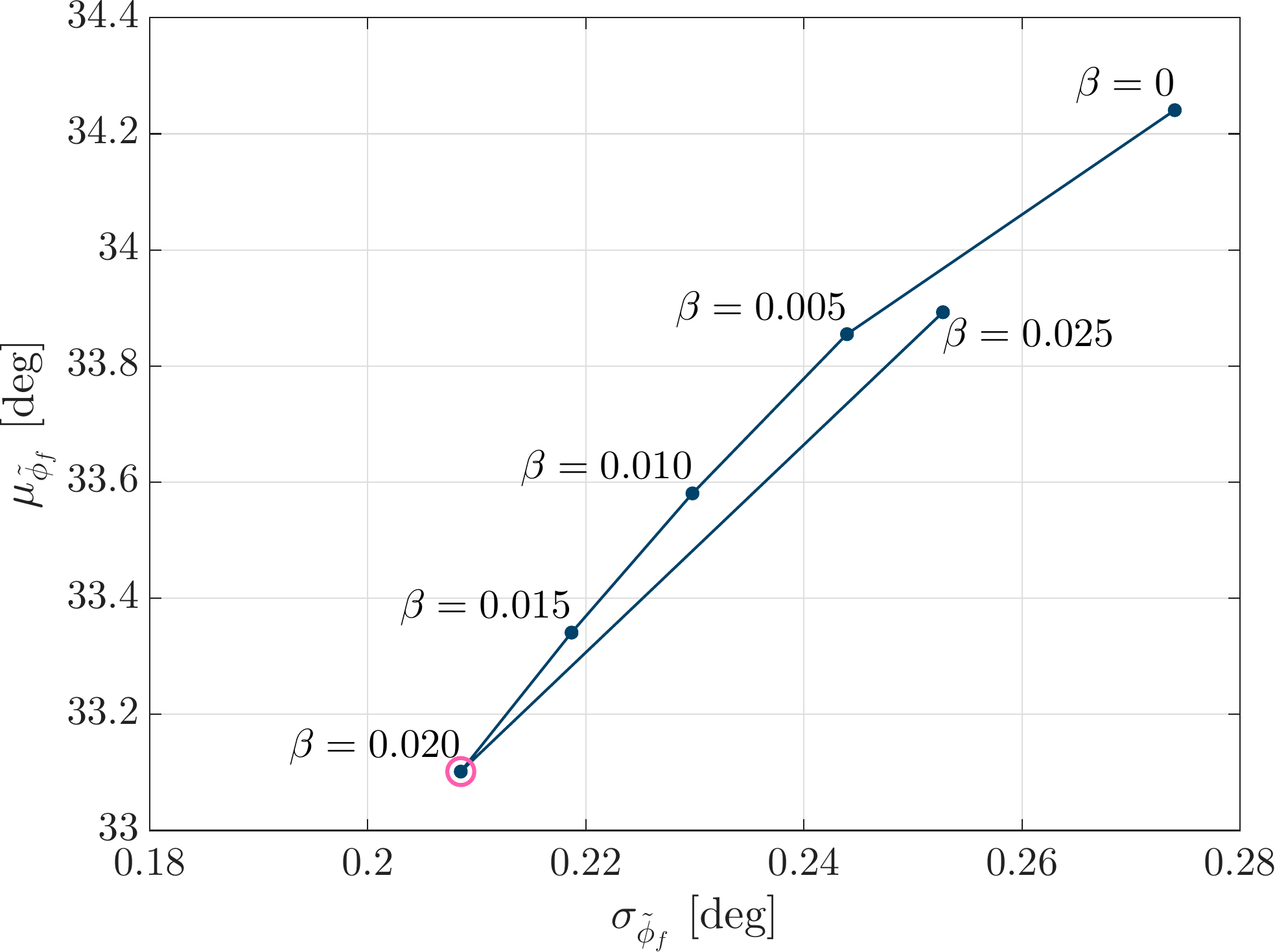}
\caption{Trade study from 100 Monte Carlo trials for $\beta \in \{0,0.005,0.010,0.015,0.020,0.025\}$.}
\label{fig: rlv-pareto-trade-study}
\end{figure}

To select the terminal desensitization weight, $\beta$, a performance vs robustness trade study was performed. Six Monte Carlo simulations were conducted corresponding to $\beta \in \{0,0.005,0.010,0.015,0.020,0.025\}$, with 100 trials per simulation. A running sensitivity penalty was also introduced with $\xi=0.5$. The mean achieved terminal latitude, ${\mu}_{\tilde{\phi}_f}$, and corresponding standard deviation, $\sigma_{\tilde{\phi}_f}$, are plotted in Fig.~\ref{fig: rlv-pareto-trade-study}. The chosen desensitization weight is located at the "knee" of the curve - the point at which no further gain in robustness can be attained. Table~\ref{tab: rlv-pareto-study} shows the mean achieved terminal latitude and corresponding standard deviation, the average computational time required for each guidance update denoted $\mu_{CPU}$, and the average percentage of successfully converged guidance updates, $\bar{R}_c$. For this study, the weight $\beta=0.02$ sacrifices 3.33\% performance to gain a 23.91\% reduction in the standard deviation. Desensitization weights larger than $\beta=0.02$ no longer reduce the standard deviation.

\begin{figure}[h!]
    \centering
    \subfloat[Altitude vs speed\label{fig: rlv-alt-vs-spd}]{
        \includegraphics[scale=0.25]{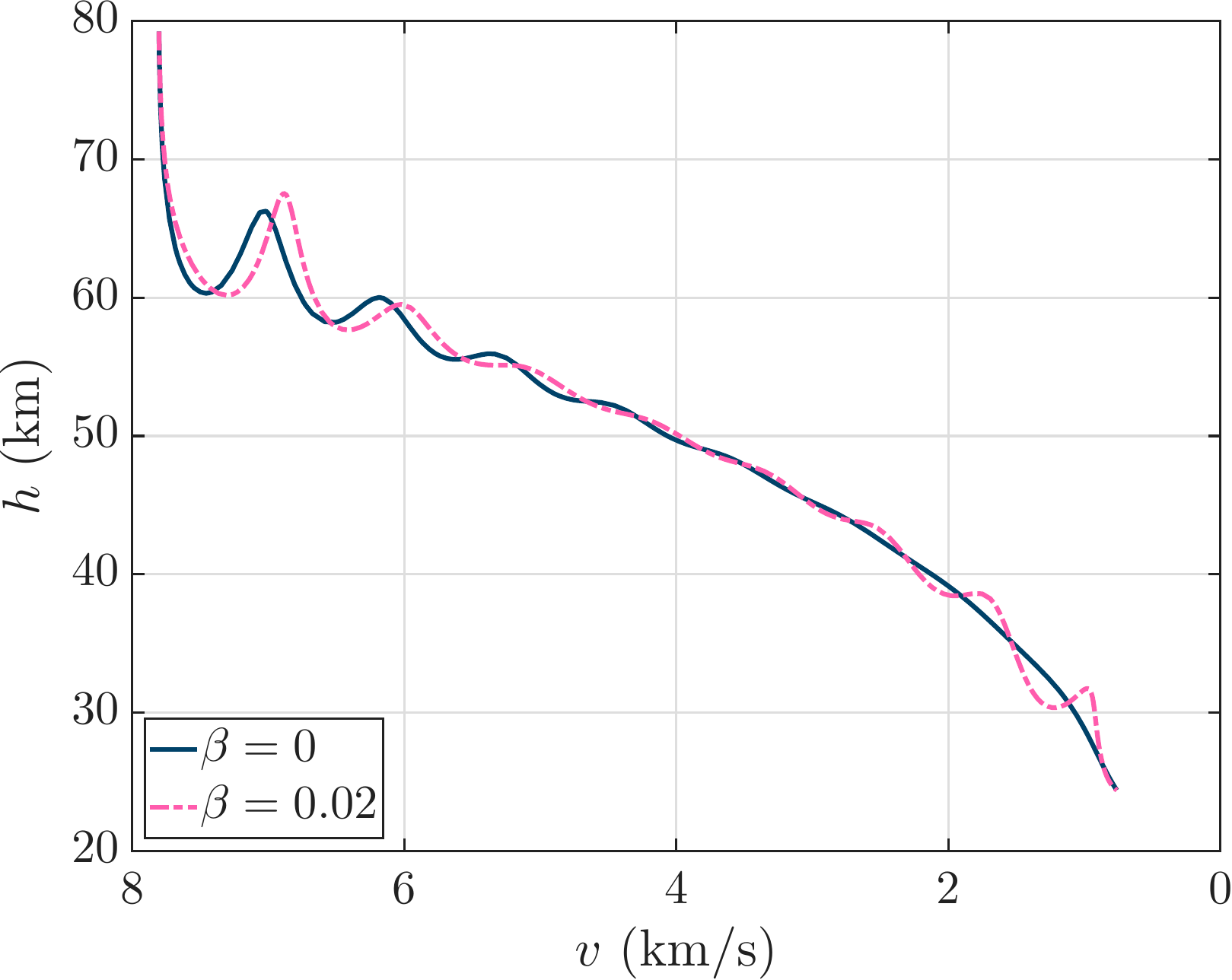}}
    \subfloat[Latitude vs longitude\label{fig: rlv-lat-vs-lon}]{
        \includegraphics[scale=0.25]{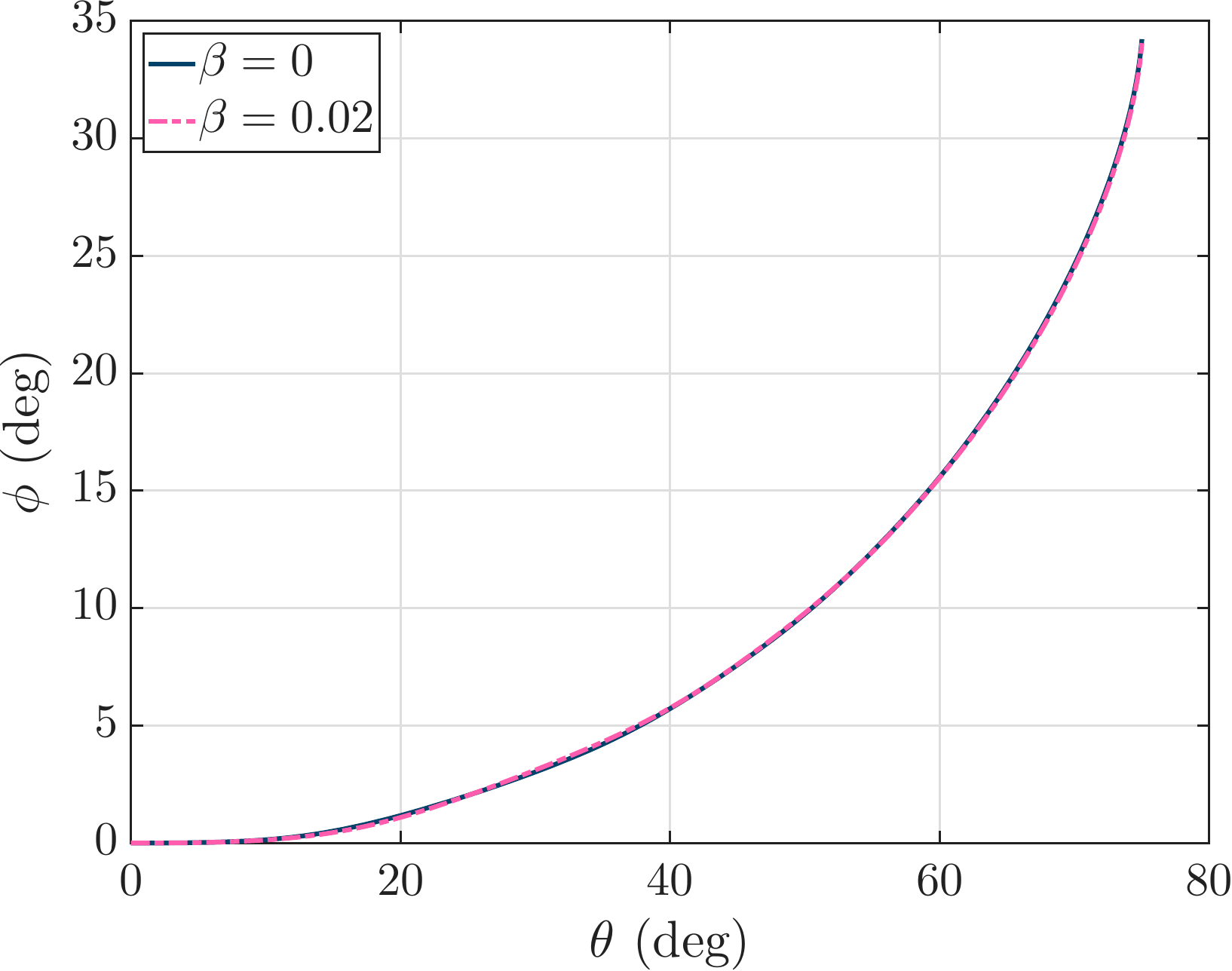}}\\
    \subfloat[Flight-path angle vs time\label{fig: rlv-fpa-vs-time}]{
        \includegraphics[scale=0.25]{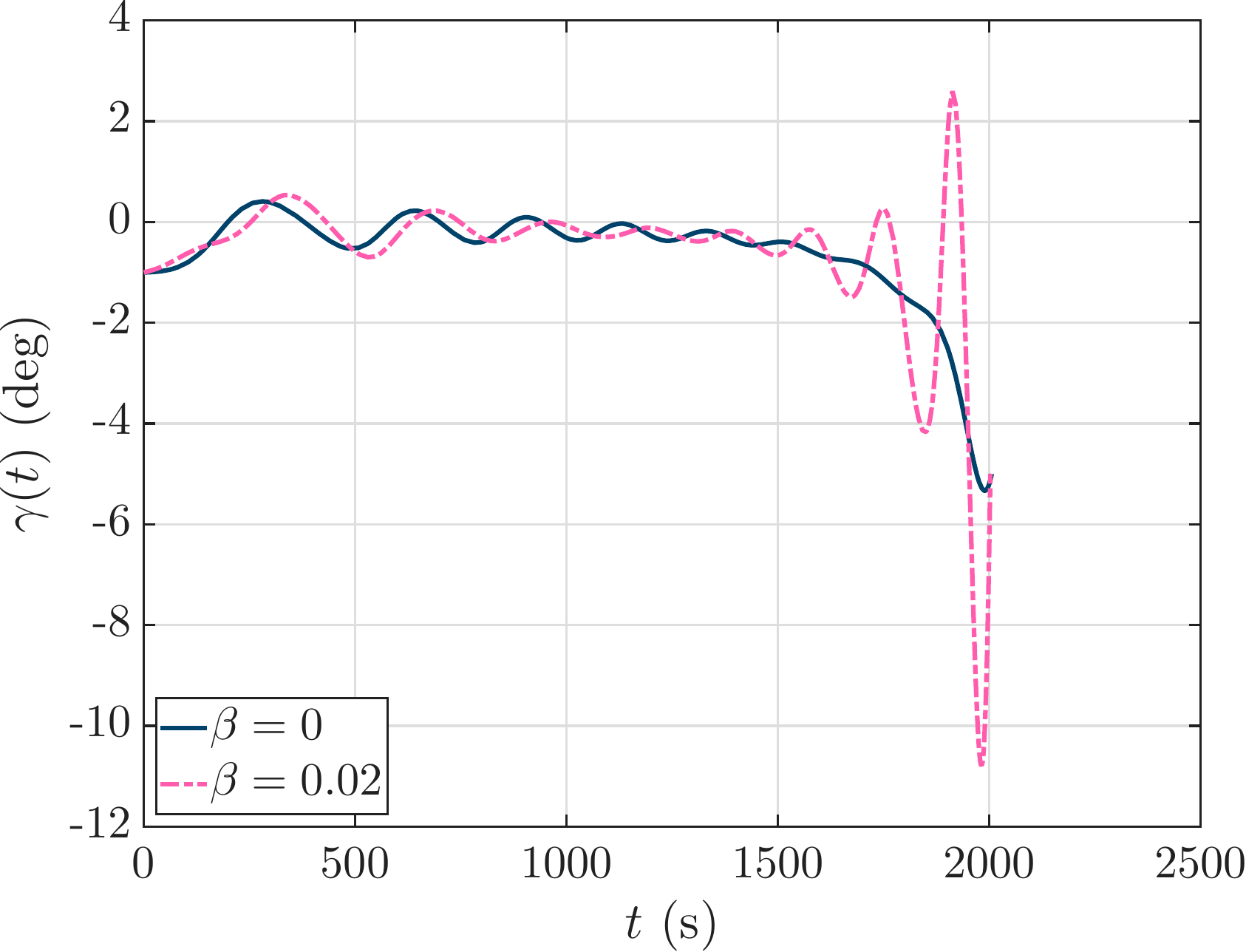}}
     \subfloat[Azimuth angle vs time\label{fig: rlv-azi-vs-time}]{
        \includegraphics[scale=0.25]{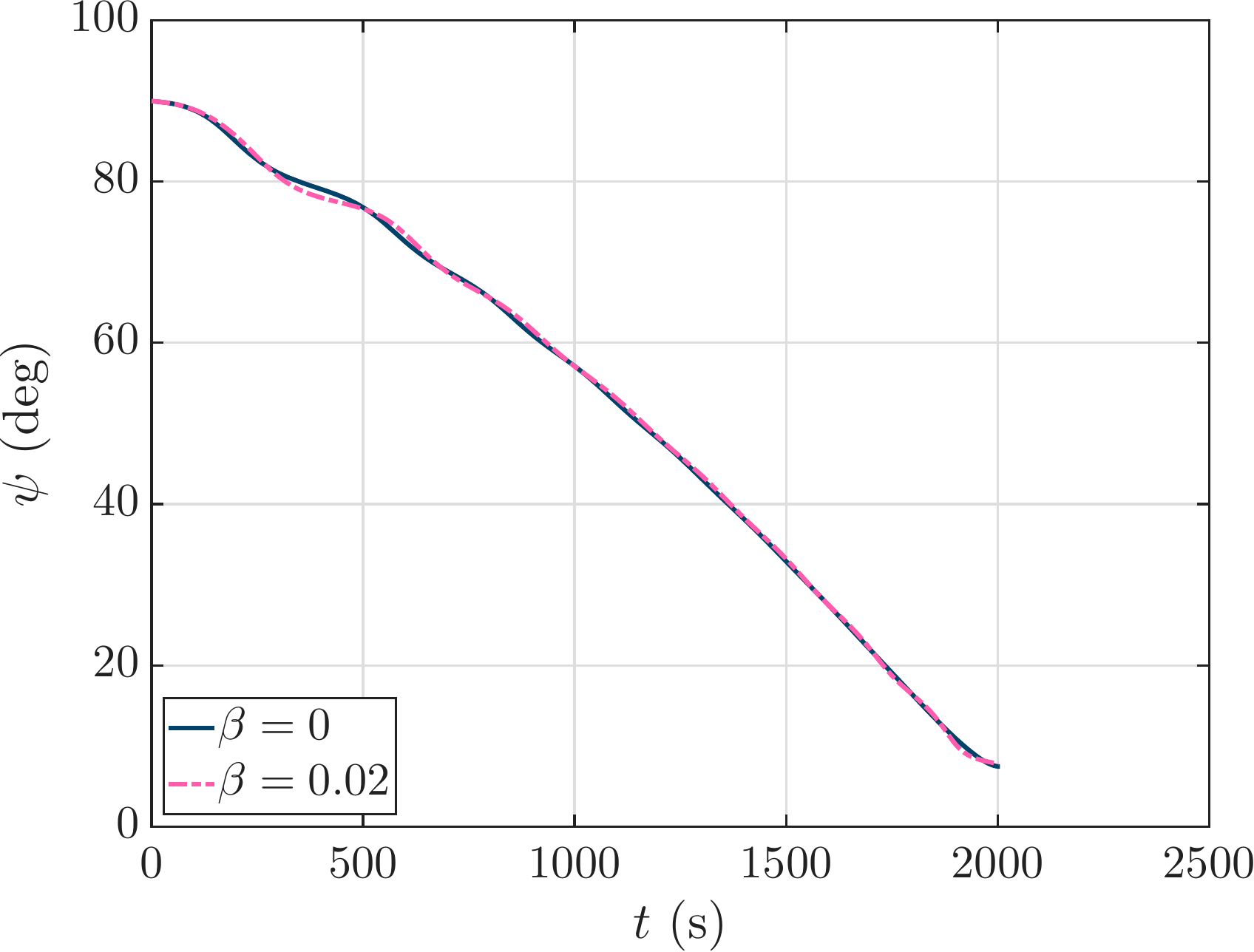}}\\
    \subfloat[Angle of attack vs time\label{fig: rlv-aoa-vs-time}]{
        \includegraphics[scale=0.25]{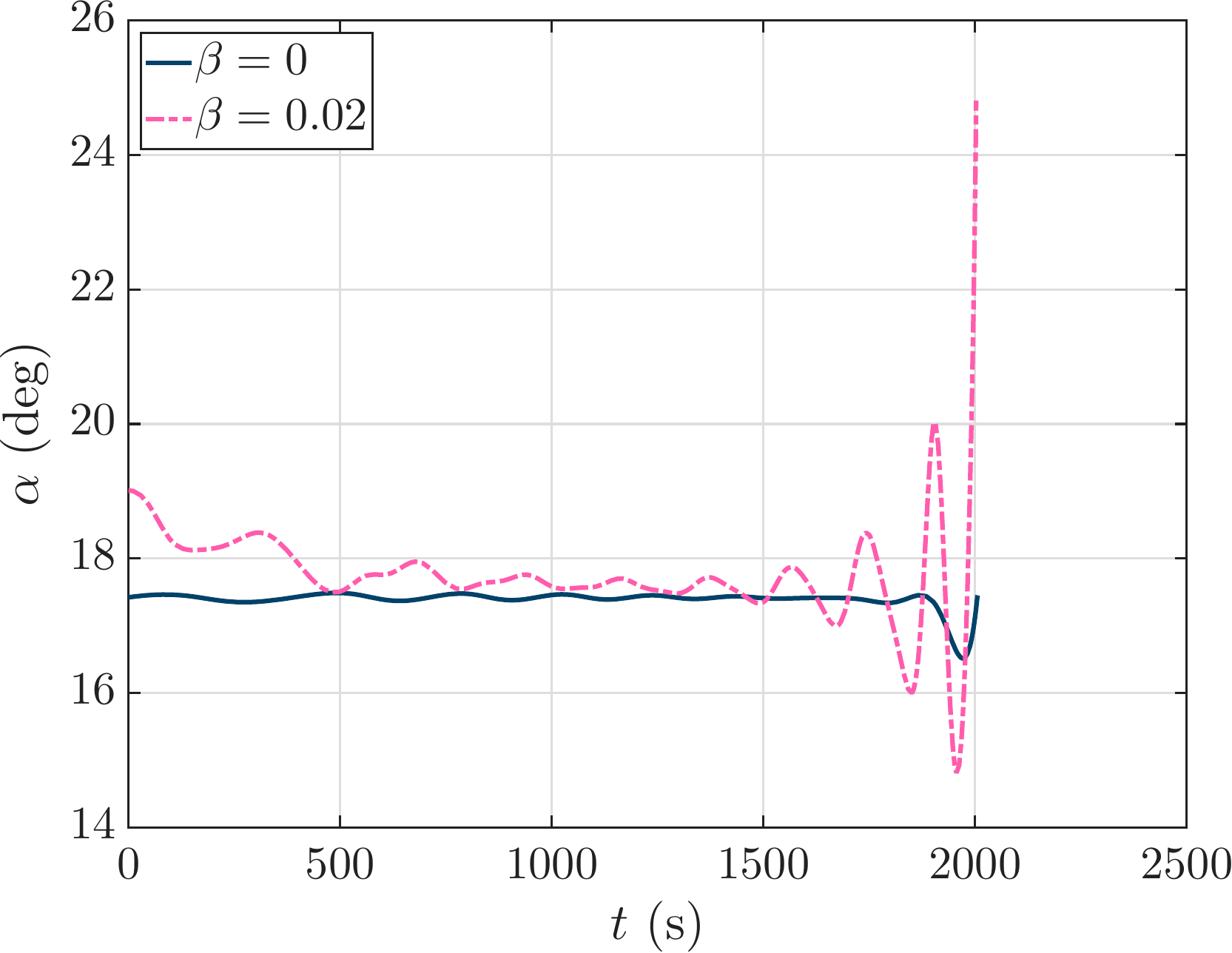}}
     \subfloat[Bank angle vs time\label{fig: rlv-bank-vs-time}]{
        \includegraphics[scale=0.25]{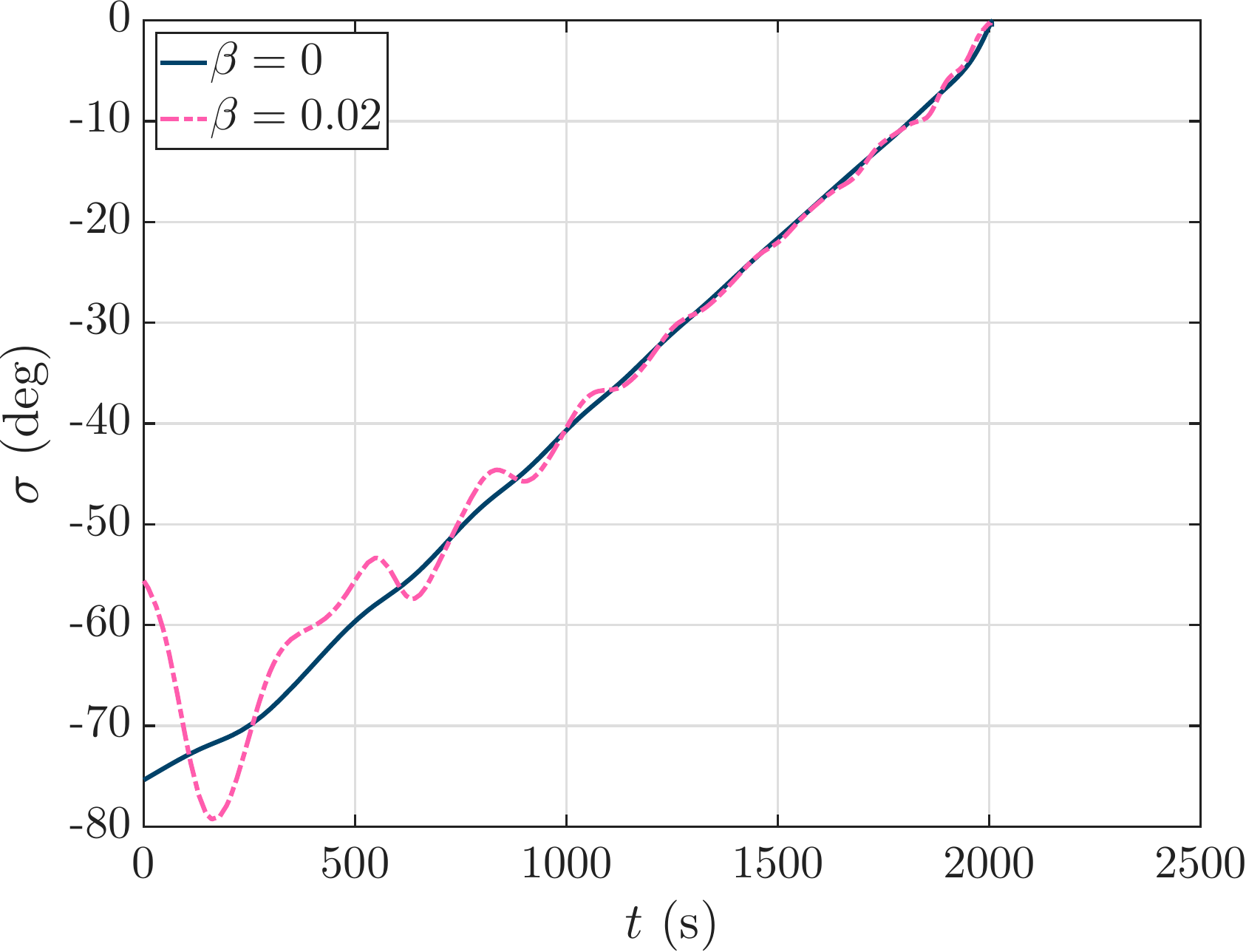}}\\
    \caption{The optimal solution for the unconstrained reusable launch vehicle entry problem with $\beta \in \{0,0.02\}$.}
    \label{fig: rlv-reference-solution}
\end{figure}

\begin{figure}[h!]
    \centering
    \subfloat[Sensitivity of $\theta$ with respect to $H$\label{fig: sens-lon-wrt-H}]{
        \includegraphics[scale=0.25]{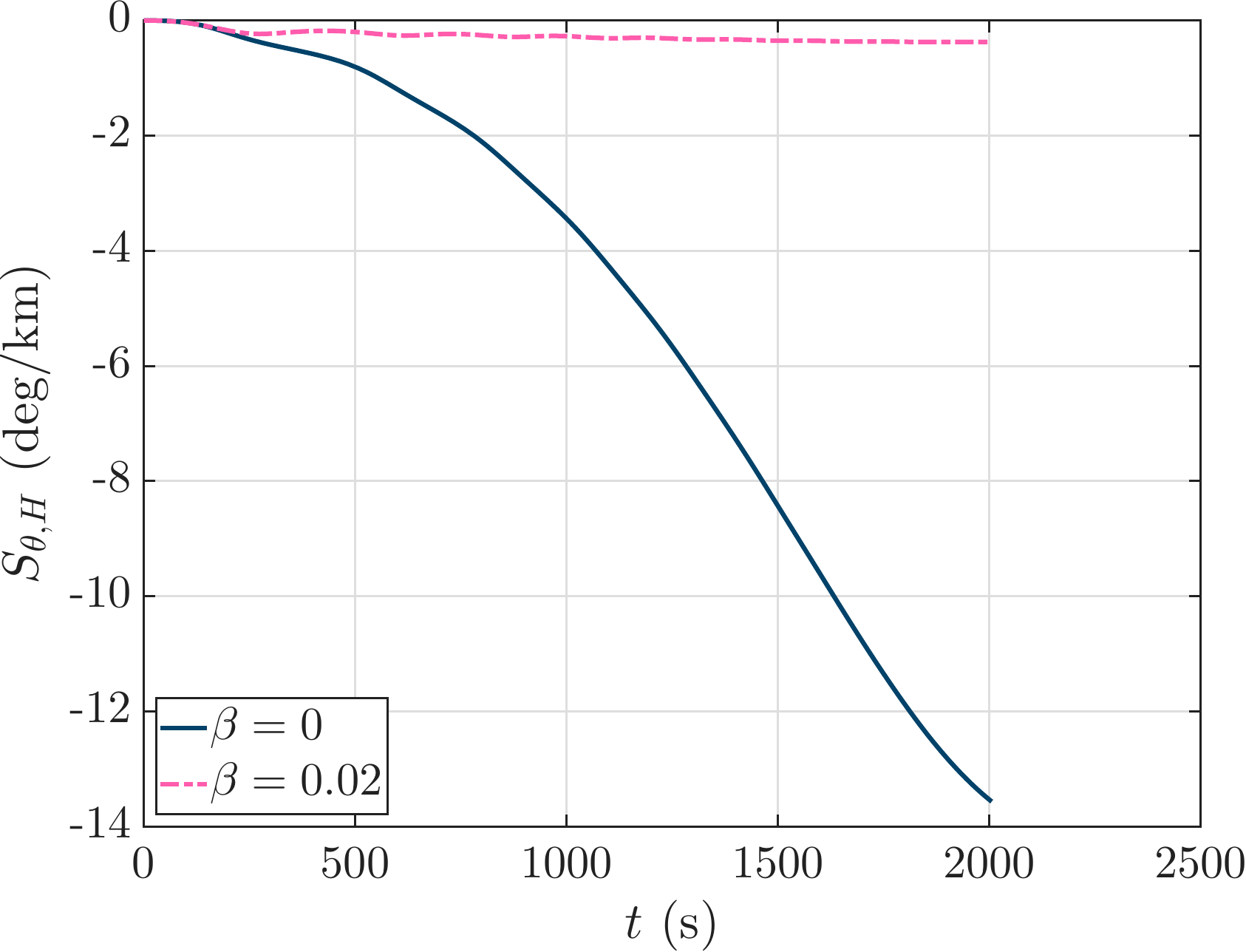}}
    \subfloat[Sensitivity of $\theta$ with respect to $C_{D_0}$\label{fig: rlv-sens-lon-wrt-cd0}]{
        \includegraphics[scale=0.25]{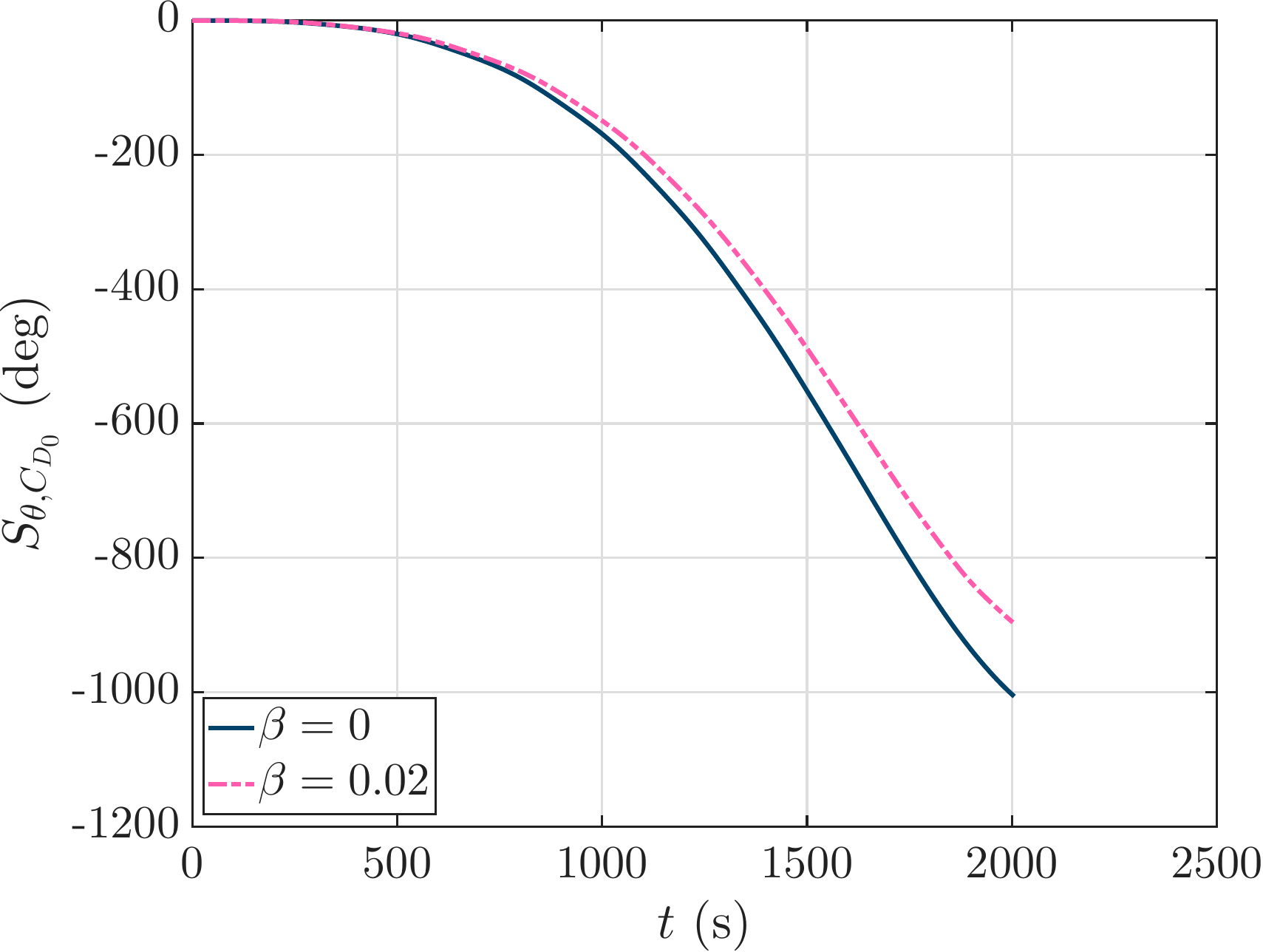}}\\
    \subfloat[Sensitivity of $\phi$ with respect to $H$\label{fig: rlv-sens-lat-wrt-H}]{
        \includegraphics[scale=0.25]{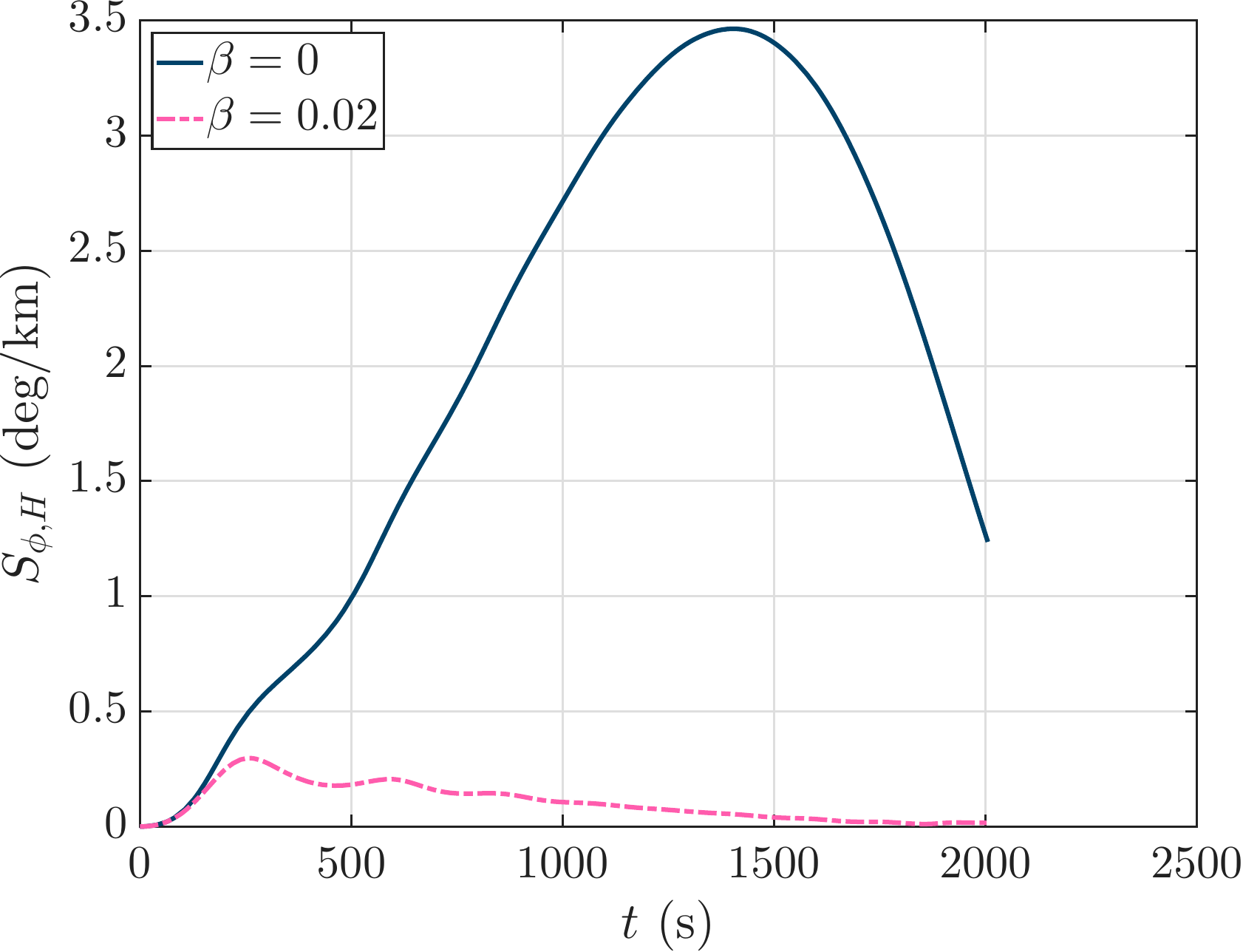}}
     \subfloat[Sensitivity of $\phi$ with respect to $C_{D_0}$\label{fig: rlv-sens-lat-wrt-cd0}]{
        \includegraphics[scale=0.25]{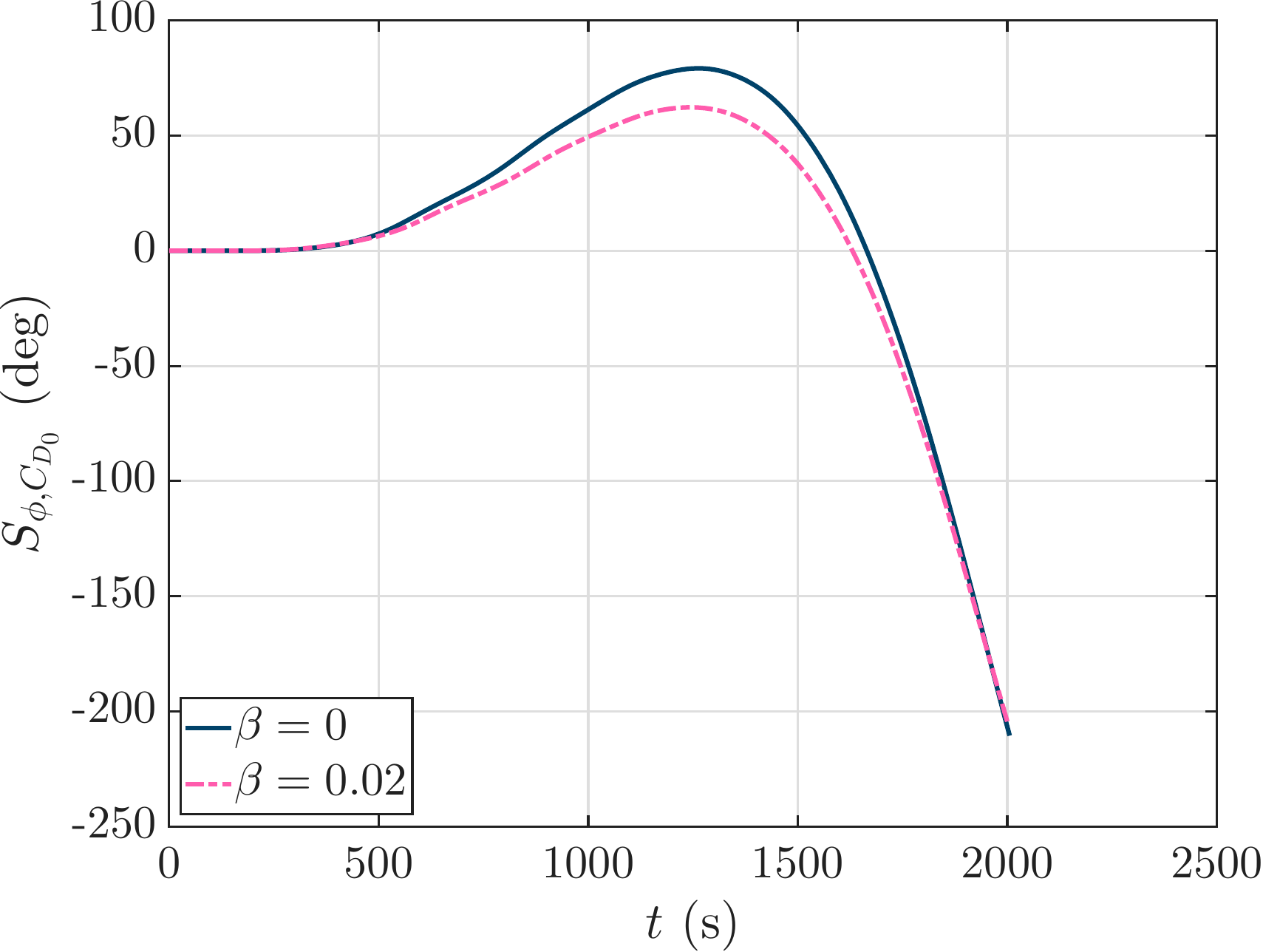}}\\
    \subfloat[Sensitivity of $v$ with respect to $H$\label{fig: rlv-sens-v-wrt-H}]{
        \includegraphics[scale=0.25]{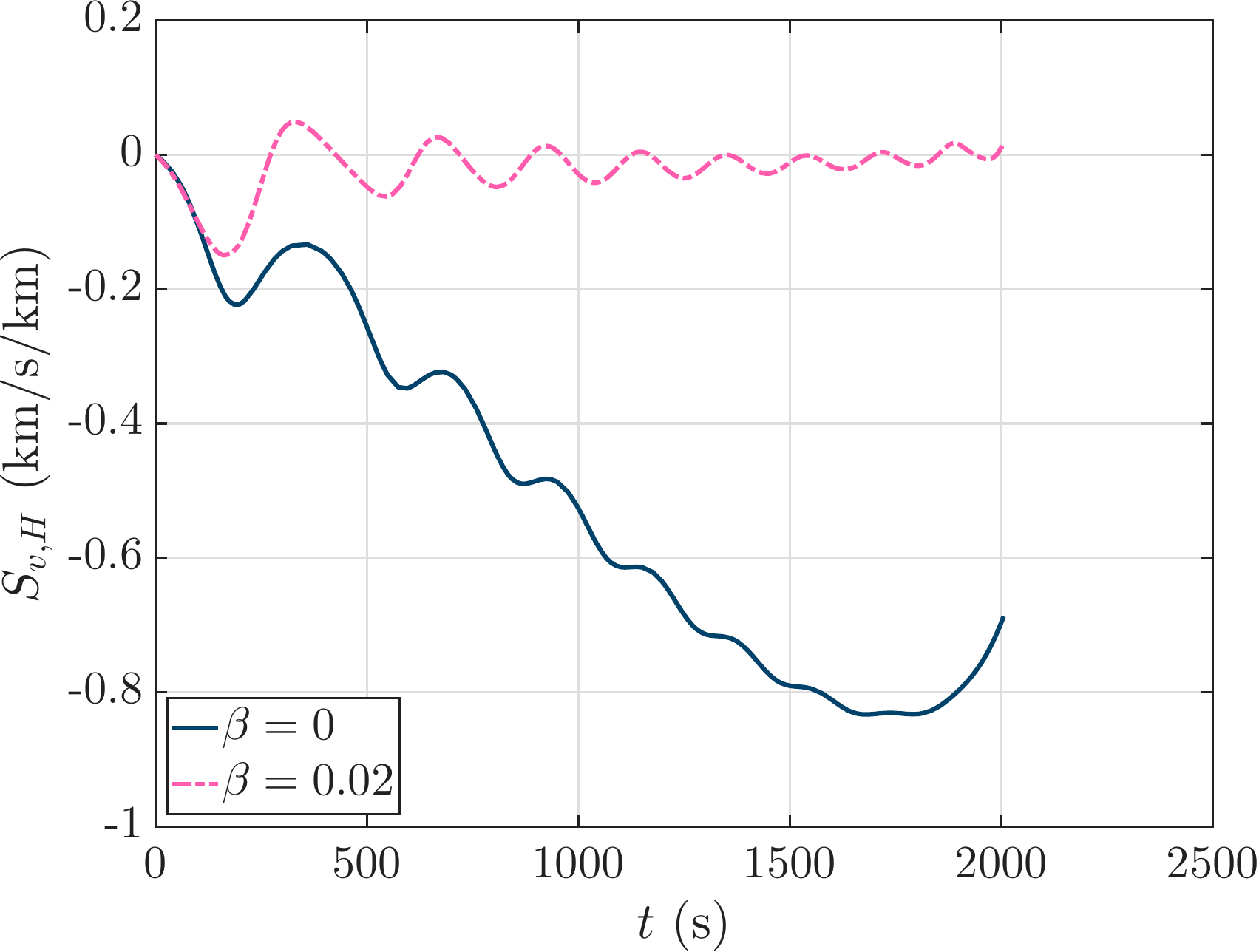}}
     \subfloat[Sensitivity of $v$ with respect to $C_{D_0}$\label{fig: rlv-sens-spd-wrt-cd0}]{
        \includegraphics[scale=0.25]{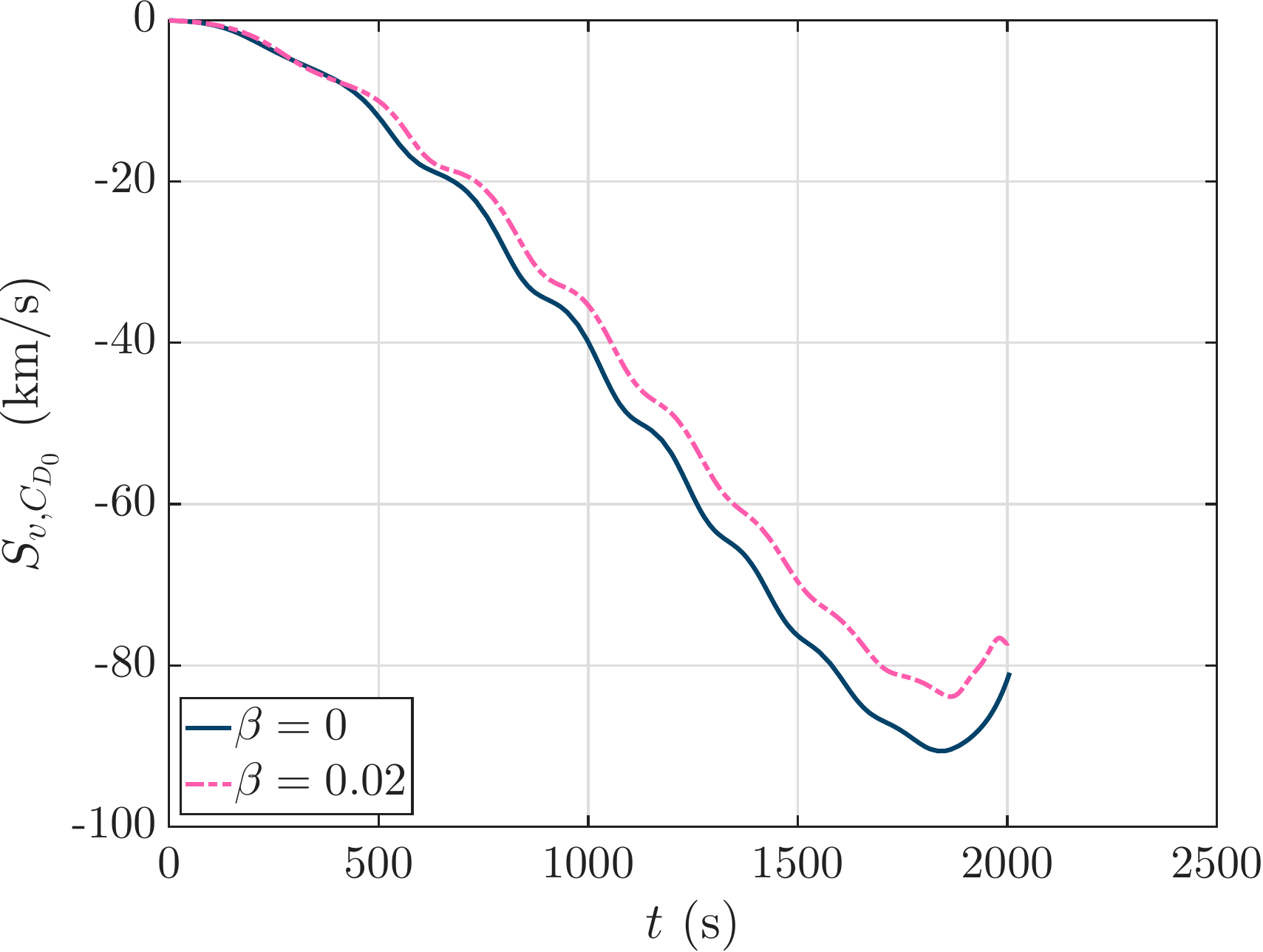}}\\
    \subfloat[Sensitivity of $\gamma$ with respect to $H$\label{fig: rlv-sens-fpa-wrt-H}]{
        \includegraphics[scale=0.25]{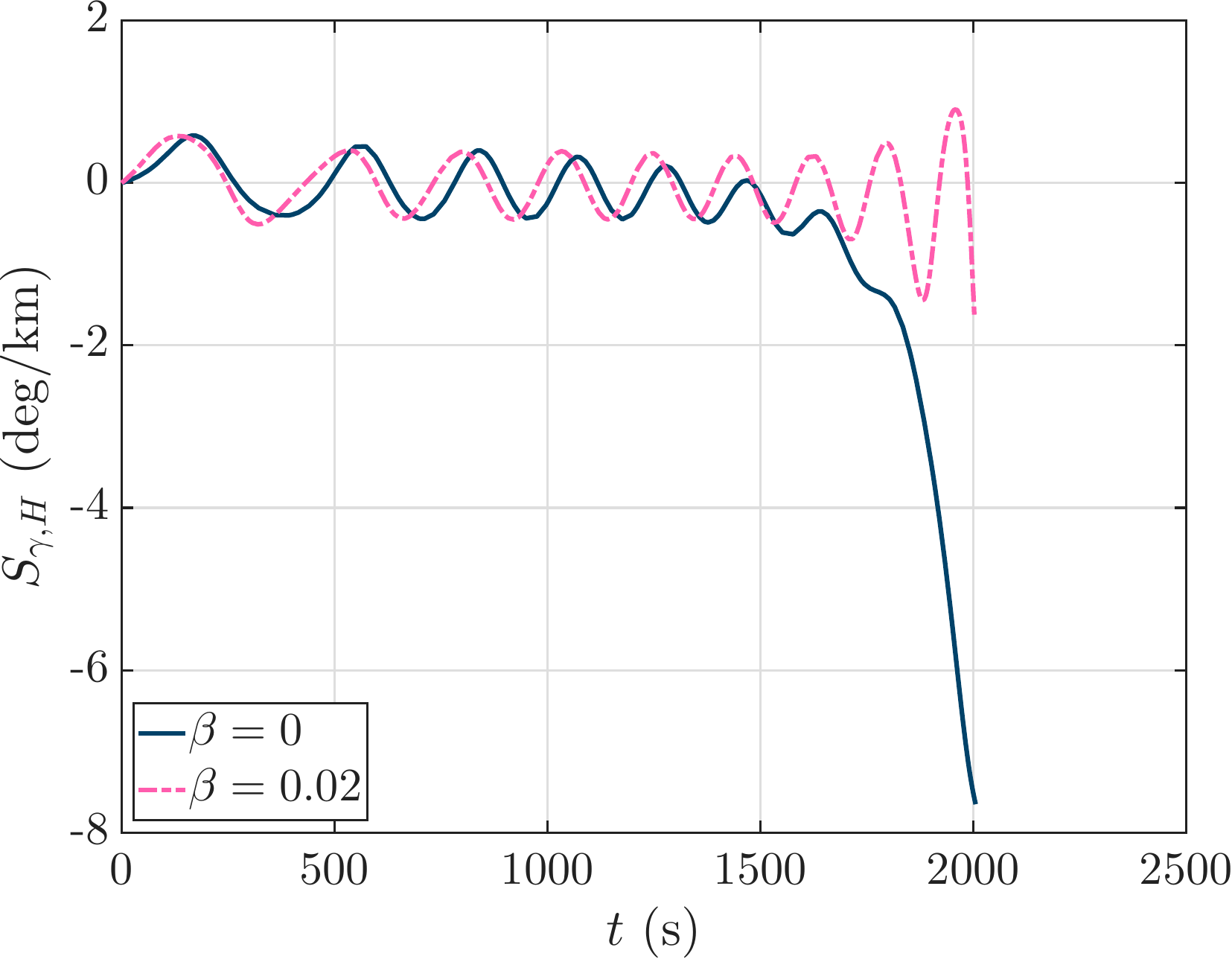}}
     \subfloat[Sensitivity of $\gamma$ with respect to $C_{D_0}$\label{fig: rlv-sens-fpa-wrt-cd0}]{
        \includegraphics[scale=0.25]{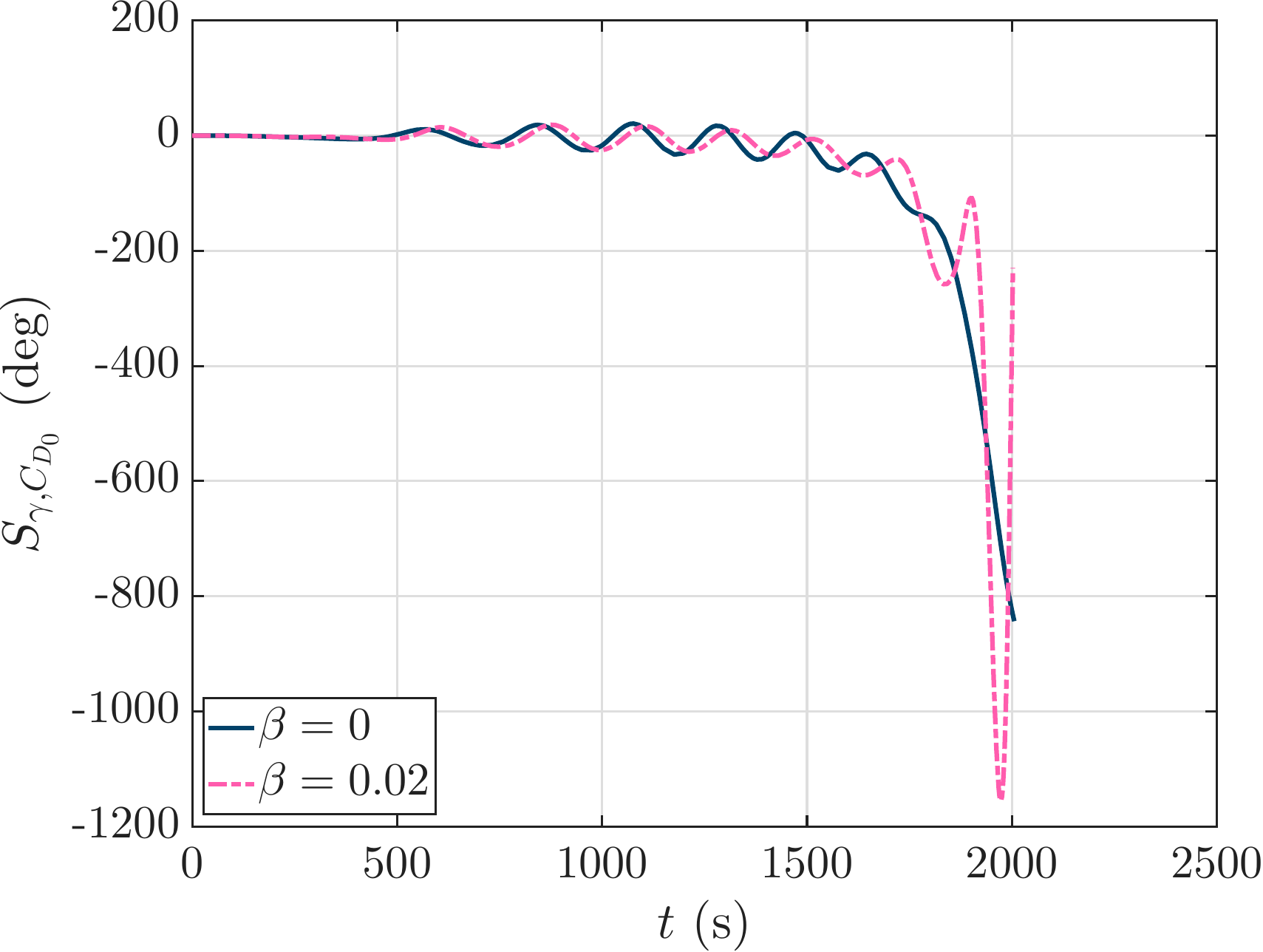}}
    \caption{Penalized state sensitivities with respect to density scale-height and zero-lift drag coefficient for $\beta \in \{0,0.02\}$.}
    \label{fig: rlv-sens-wrt-p}
\end{figure}

\begin{figure}[h!]
    \centering
    \subfloat[Sampled parameter 99.97\% confidence ellipse\label{fig: rlv-parameters}]{
        \includegraphics[scale=0.36]{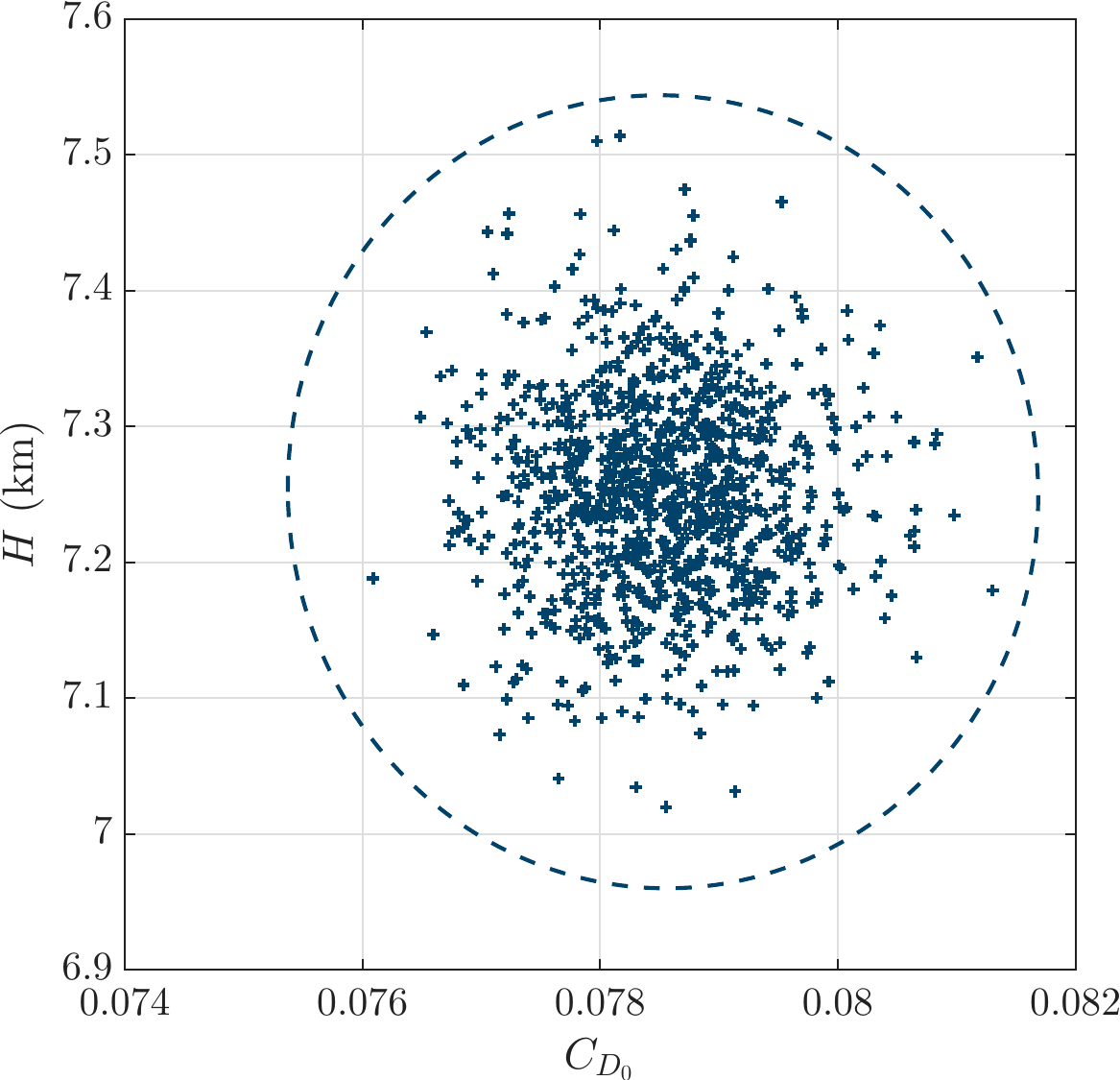}}
    \subfloat[Crossrange-downrange 99.97\% confidence ellipse \label{fig: rlv-terminal-deviation}]{
        \includegraphics[scale=0.36]{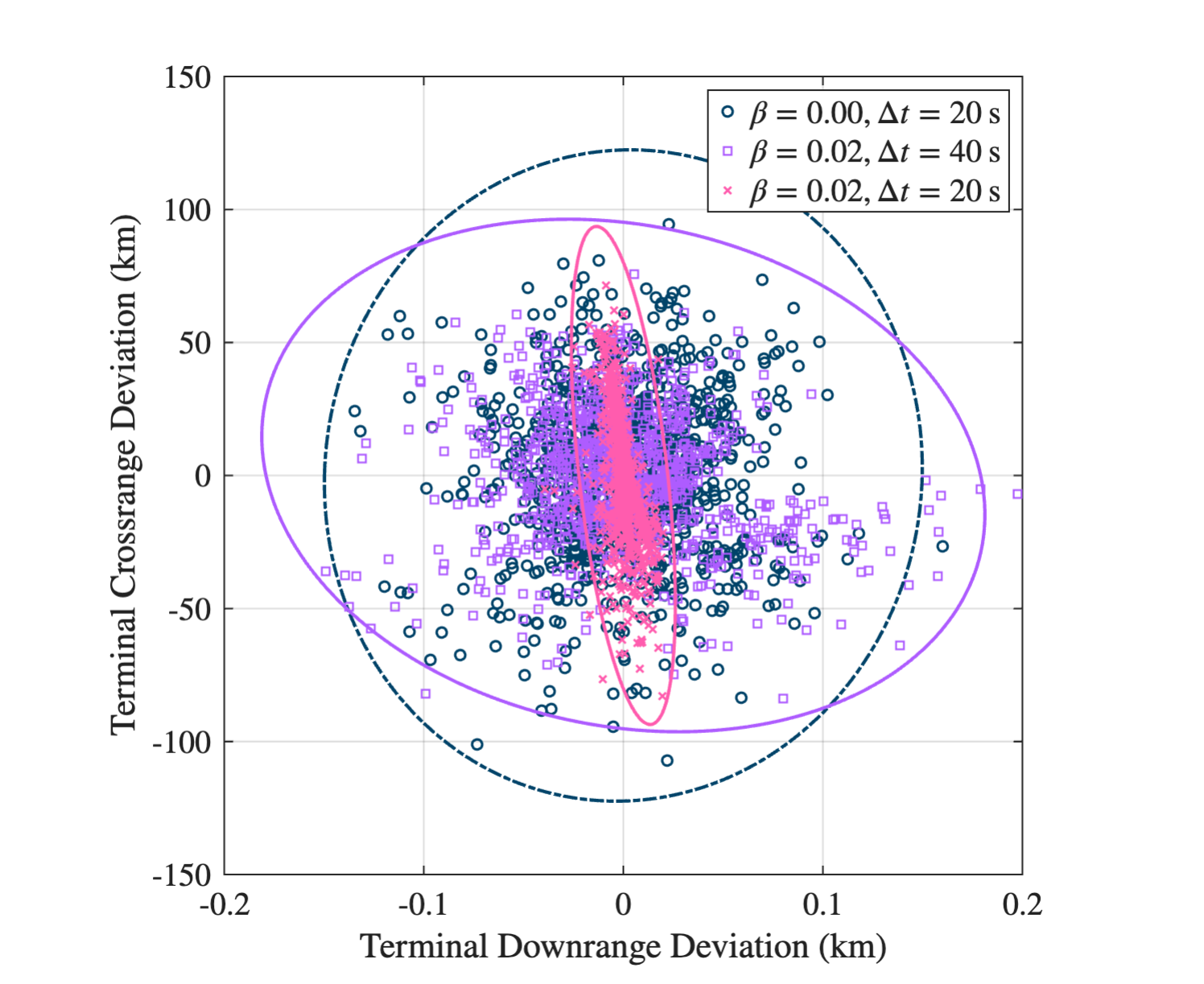}}\\
    \subfloat[Crossrange-downrange trajectory envelope \label{fig: rlv-crossrange-vs-downrange}]{
        \includegraphics[scale=0.34]{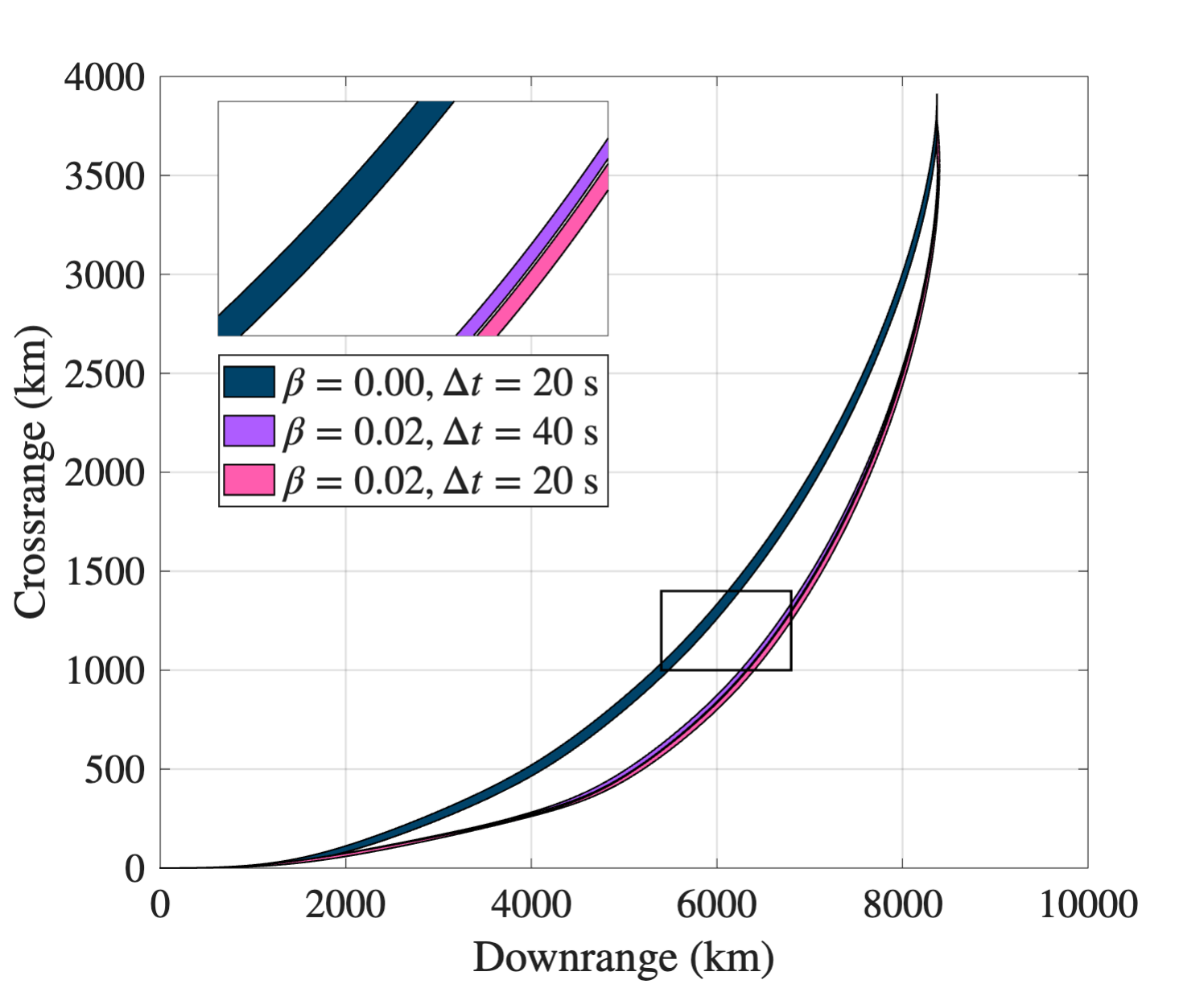}}
    \caption{Results from 1000 Monte Carlo trials for terminal desensitization weight $\beta \in \{0,0.02\}$ and guidance cycle duration $\Delta t \in \{20,40\}$ seconds.}
    \label{fig: rlv-MC-results}
\end{figure}

\begin{table}[h!]
\caption{\label{tab: rlv-pareto-study} Trade study from 100 Monte Carlo trials for $\beta \in \{0,0.005,0.010,0.015,0.020,0.025\}$.}
\centering
\begin{tabular}{cccccc}
\hline
$\beta$ & $\mu_{\tilde{\phi}_f}$ (deg) & $\sigma_{\tilde{\phi}_f}$ (deg) & $\mu_{CPU}$ (s) & $\bar{R}_c$ (\%)\\
\hline
0.000 & 34.2407 & 0.2740 & 0.6010 & 100.0\\
0.005 & 33.8549 & 0.2439 & 1.0774 & 100.0\\
0.010 & 33.5806 & 0.2298 & 1.3342 & 100.0\\
0.015 & 33.3409 & 0.2187 & 1.9101 & 100.0\\
0.020 & 33.1015 & 0.2085 & 1.7043 & 99.99\\
0.025 & 33.8927 & 0.2527 & 2.1997 & 100.0\\
\hline
\end{tabular}
\end{table}
For a more in-depth analysis of the proposed DOG framework, three additional Monte Carlo simulations were performed with 1000 trials per simulation. The first simulation is representative of no desensitization ($\beta=0$), and the second simulation corresponds to the "tuned" desensitization weight, $\beta=0.02$. Both the first and second simulations perform guidance updates at 20-second intervals. The third simulation, however, employs a desensitization weight of $\beta=0.02$ at 40-second intervals. The purpose of the third simulation is to demonstrate the computational efficiency of DOG when fewer control corrections are allowed. The reference optimal state and control are presented in Fig.~\ref{fig: rlv-reference-solution}. Figure~\ref{fig: rlv-sens-wrt-p} shows the sensitivities of the penalized state variables with respect to the density scale-height and zero-lift drag coefficient. The perturbed parameter values along with the terminal state error distribution for the penalized crossrange and downrange are shown in Fig.~\ref{fig: rlv-parameters} and ~\ref{fig: rlv-terminal-deviation}, respectively. Figure~\ref{fig: rlv-crossrange-vs-downrange} features the Monte Carlo trajectory envelopes for the penalized crossrange-downrange trajectories when guidance updates are permitted at 20-second and 40-second intervals. Table~\ref{tab: rlv-refined-pareto-study} documents the achieved mean performance index for each simulation along with the standard deviation of the terminal latitude and longitude, the mean CPU time, and the mean convergence rate.

\begin{table}[htb!]
\caption{\label{tab: rlv-refined-pareto-study} Results from 1000 Monte Carlo trials for $\beta \in \{0,0.02\}$.}
\centering
\begin{tabular}{cccccccc}
\hline
$\beta$ & $\Delta t$ (s) & $\mu_{\tilde{\phi}_f}$ (deg) & $\sigma_{\tilde{\phi}_f}$ (deg) & $\sigma_{\tilde{\theta}_f}$ (deg) & $\mu_{CPU}$ (s) & $\bar{R}_c$ (\%)\\
\hline
0.00 & 20 & 34.2099 & 0.2723 & 0.4030 & 0.6158 & 100\\
0.02 & 20 & 33.0780 & 0.2082 & 0.0691 & 1.8026 & 100\\
0.02 & 40 & 33.2490 & 0.2143 & 0.4819 & 2.2270 & 99.99\\
\hline
\end{tabular}
\end{table}

\clearpage
\subsection{Discussion of Results}\label{section: results-discussion}
The two examples discussed in Section~\ref{section: numerical-examples} demonstrate the key features of the proposed desensitized optimal guidance framework. The principal feature of this method is a reduction in the dispersion of penalized state variables in the presence of parametric uncertainty. By augmenting the objective functional to include sensitivity penalties, computed reference solutions are inherently robust to model uncertainties. As the weight associated with this sensitivity penalty increases, the control solution trades performance for robustness. 
This trade is most clearly illustrated in Zermelo's navigation problem. For the largest desensitization weight evaluated, $\beta=0.40$, the desensitized solution sacrifices 17.62\% of the terminal horizontal distance, $x_1(t_f)$, in exchange for a 36.80\% reduction in the corresponding standard deviation. The RLV Entry problem of Section~\ref{section: rlv-entry} exhibits a similar behavior. For $\beta=0.02$, a 3.31\% reduction in the terminal latitude yields a 23.54\% reduction in the associated standard deviation. This relatively large growth in robustness is partly attributed to the presence of both a \textit{running} and \textit{terminal} sensitivity penalty within the augmented cost. An additional measure of performance is the terminal distribution for the penalized state variables that are \textit{fixed} at the final time. In the RLV entry problem, the coupled longitude and latitude dynamics are largely influenced by the aerodynamic and atmospheric model uncertainties. Thus, the benefits of DOG are most clearly evident in the terminal longitude distribution. For $\beta=0.02$, Monte Carlo analysis revealed an 82.85\% reduction in the standard deviation of the achieved terminal longitude. Both examples demonstrate that the proposed framework significantly improves the precision of the realized terminal states with minimal degradation in performance.

The degree to which DOG can mitigate trajectory variations is governed by the magnitude of the desensitization weight. The desensitization weight can be selected such that a balanced trade-off between performance and robustness can be realized, where tools such as Pareto analysis can be leveraged to inform an appropriate value. For the RLV entry problem, the Pareto study determines the desensitization weight that minimizes the standard deviation of the terminal latitude while maintaining acceptable nominal performance. This optimal balance is often found at the "knee" of the Pareto curve where further increases in the desensitization weight yield diminishing reductions in dispersion, as indicated in Fig.~\ref{fig: rlv-pareto-trade-study}. As this weight increases, desensitization takes more precedence in the OCP, as illustrated in Fig.~\ref{fig: rlv-sens-wrt-p}. Consequently, the resulting trajectories exhibit less sensitivity to model imperfections, leading to a tighter trajectory envelope observed through the Monte Carlo simulation in Figs.~\ref{fig: zermelo-envelope} and \ref{fig: rlv-crossrange-vs-downrange}. In Fig.~\ref{fig: rlv-crossrange-vs-downrange}, enforcing a desensitization weight of $\beta=0.02$ in conjunction with guidance corrections at 20-second intervals tightens the maximum width of the trajectory envelope by approximately 37 km (45.89\% reduction) while the average guidance update CPU time remains less than two seconds. 

In addition to improving robustness to parametric uncertainty, the proposed framework enables guidance updates to be performed less frequently without significantly degrading trajectory precision. Figure~\ref{fig: rlv-terminal-deviation} compares the terminal crossrange-downrange dispersion for the nominal solution ($\beta=0$) with updates performed in 20-second intervals and the desensitized solution ($\beta=0.02$) with updates performed in 40-second intervals. 
The case where $\beta=0$ corresponds to the computational optimal guidance approach of Ref.~\cite{dennis-hager-2019}, where the optimal trajectory is repeatedly recomputed without sensitivity augmentation.
Despite the guidance update frequency being roughly halved, the desensitized solution still achieved a 21.30\% reduction in the standard deviation of the terminal latitude in exchange for a 2.81\% loss in nominal performance. Although the standard deviation in the terminal longitude increased by 19.56\% for the desensitized solution, the crossrange-downrange trajectory envelope remained tighter (approximately 40 km less) than the Monte Carlo envelope produced by the nominal solution with guidance updates every 20 seconds. Collectively, the results of this study demonstrate that DOG can improve system robustness to model uncertainties while reducing the required guidance update frequency, supporting the feasibility of DOG for real-time onboard implementation.

The method depends upon the weights that are chosen, the guidance cycle time, and the size of the disturbance. First, it is noted that the sensitivity weights are problem dependent and must be chosen appropriately for a given application. Second, the duration of the guidance cycle is inversely proportional to the magnitudes of the disturbances. When the disturbances are known to be smaller, the guidance cycle duration can be made larger. Conversely, when the disturbances are larger, the guidance cycle duration must made smaller. All of these choices are design parameters that are specific to the application under consideration.

\section{Conclusions}\label{section: conclusion}
A computational method for the desensitized optimal guidance of dynamical systems subject to parametric uncertainty has been developed. The method directly integrates a reduced desensitized trajectory optimization technique within the guidance update process to mitigate trajectory dispersion. A previously developed mesh truncation and mesh realignment strategy is used to enable rapid convergence to the optimal solution. The method has been demonstrated on two numerical examples. Results show that explicitly penalizing state sensitivities within the objective functional reduces trajectory dispersion of the penalized state variables \textit{and} variability of the realized performance index with modest degradation of the nominal cost. Through Monte Carlo analysis, it is found that desensitized guidance solutions allow guidance updates to be performed less frequently without significantly degrading trajectory precision, contributing to the computational feasibility of the proposed guidance framework for onboard implementation. Together, these results highlight the potential of desensitized optimal guidance for robust real-time trajectory generation in the presence of model uncertainties.

\section*{Acknowledgments}
The authors acknowledge support for this research from the U.S.
Office of Naval Research under grant N00014-22-1-2397, from the
U.S. National Science Foundation under grant CMMI-2031213, and from the U.S. Air Force Research Laboratory under grant FA8651-24-1-0004.

\end{document}